# ELECTRODYNAMICS and the GAUSS LINKING INTEGRAL

## on the 3-SPHERE and in HYPERBOLIC 3-SPACE

*Dennis DeTurck and Herman Gluck*

In this first of two papers, we develop a steady-state version of classical electrodynamics on the 3-sphere and in hyperbolic 3-space, including an explicit formula for the vector-valued Green's operator, an explicit formula of Biot-Savart type for the magnetic field, and a corresponding Ampere's Law contained in Maxwell's equations.

We then use this to obtain explicit integral formulas for the linking number of two disjoint closed curves in these spaces.

All the formulas, like their prototypes in Euclidean 3-space, are geometric rather than just topological, in the sense that their integrands are invariant under orientation-preserving isometries of the ambient space.

In the second paper, we obtain integral formulas for twisting, writhing and helicity, and prove the theorem LINK = TWIST + WRITHE in the 3-sphere and in hyperbolic 3-space. We then use these results to derive upper bounds for the helicity of vector fields and lower bounds for the first eigenvalue of the curl operator on subdomains of these two spaces.

An announcement of most of these results, and a hint of their proofs, can be found in the Math ArXiv, math.GT/0406276.

# ORGANIZATION

The integral formulas in this paper contain vectors lying in different tangent spaces; in non-Euclidean settings these vectors must be moved to a common location to be combined.

In $S^3$ regarded as the group of unit quaternions, equivalently as $SU(2)$, the differential $L_{yx^{-1}}$ of left translation by $yx^{-1}$ moves tangent vectors from $x$ to $y$. In either $S^3$ or $H^3$, parallel transport $P_{yx}$ along the geodesic segment from $x$ to $y$ also does this. As a result, we get three versions for each of the formulas that appear in the theorems below.

In part A of this paper, we present the integral formulas for linking, for the steady-state magnetic field, and for the vector-valued Green's operator in $S^3$ and $H^3$, as well as the statement of Maxwell's equations in these two spaces.

In part B, we prove the formulas for the magnetic field and the Green's operator on $S^3$ in left translation format.

In part C, we prove these formulas on $S^3$ in parallel transport format.

In part D, we prove them in $H^3$ in parallel transport format.

In part E, we verify Maxwell's equations and then prove the linking formulas in $S^3$ and $H^3$ in all formats.

# A. STATEMENTS OF RESULTS

## 1. Linking integrals in $R^3$, $S^3$ and $H^3$.

Let $K_1 = \{x(s)\}$ and $K_2 = \{y(t)\}$ be disjoint oriented smooth closed curves in either Euclidean 3-space $R^3$, the unit 3-sphere $S^3$ or hyperbolic 3-space $H^3$, and let $\alpha(x, y)$ denote the distance from $x$ to $y$.

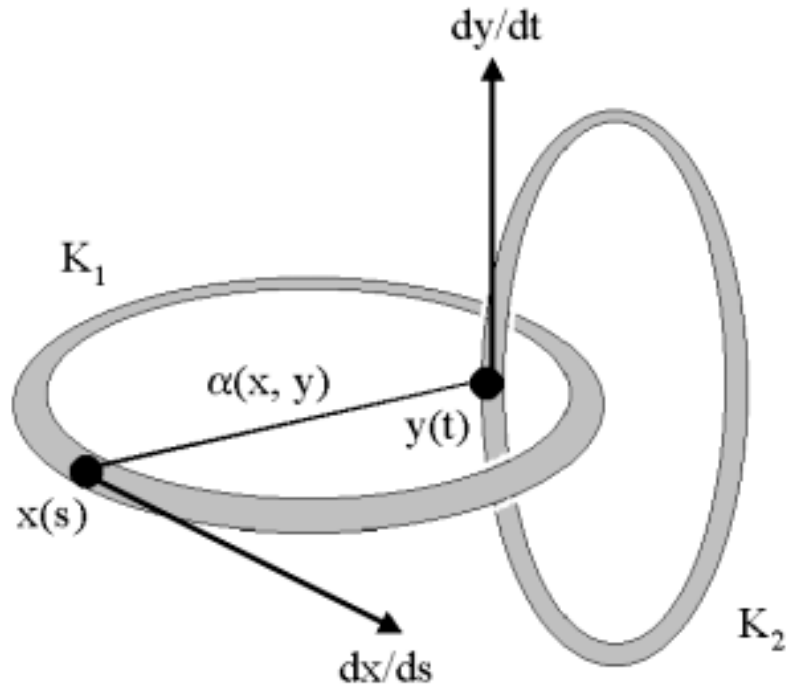

**Two linked curves**

The linking number of these two curves is defined to be the intersection number of either one of them with an oriented surface bounded by the other. It is understood that the ambient space is also oriented. The linking number depends neither on the choice of surface, nor on which of the curves is used to bound the surface.

Carl Friedrich Gauss, in a half-page paper dated January 22, 1833, gave an integral formula for the linking number in Euclidean 3-space,

$$Lk(K_1, K_2) = 1/(4\pi) \int_{K_1 \times K_2} dx/ds \times dy/dt \bullet (x - y)/|x - y|^3 \; ds \, dt \, .$$

It will be convenient for us to write this as

$$Lk(K_1, K_2) = 1/(4\pi) \int_{K_1 \times K_2} dx/ds \times dy/dt \bullet \nabla_y \varphi(x, y) \; ds \, dt \, ,$$

where $\varphi(\alpha) = 1/\alpha$, and where we use $\varphi(x, y)$ as an abbreviation for $\varphi(\alpha(x, y))$. The subscript $y$ in the expression $\nabla_y \varphi(x, y)$ tells us that the differentiation is with respect to the $y$ variable.

The following theorem gives the corresponding linking integrals on the 3-sphere and in hyperbolic 3-space. Since the location of tangent vectors is now important, we note that the vector $\nabla_y \varphi(x, y)$ is located at the point $y$.

**THEOREM 1. LINKING INTEGRALS in $S^3$ and $H^3$.**

**(1)** *On $S^3$ in left-translation format:*

$$Lk(K_1, K_2) = 1/(4\pi^2) \int_{K_1 \times K_2} L_{yx^{-1}} dx/ds \times dy/dt \bullet \nabla_y \varphi(x, y) \; ds \, dt$$
$$- 1/(4\pi^2) \int_{K_1 \times K_2} L_{yx^{-1}} dx/ds \bullet dy/dt \; ds \, dt \, ,$$

*where* $\varphi(\alpha) = (\pi - \alpha) \cot \alpha$.

**(2)** *On $S^3$ in parallel transport format:*

$$Lk(K_1, K_2) = 1/(4\pi^2) \int_{K_1 \times K_2} P_{yx} dx/ds \times dy/dt \bullet \nabla_y \varphi(x, y) \; ds \, dt \, ,$$

*where* $\varphi(\alpha) = (\pi - \alpha) \csc \alpha$.

**(3)** *On $H^3$ in parallel transport format:*

$$Lk(K_1, K_2) = 1/(4\pi) \int_{K_1 \times K_2} P_{yx} dx/ds \times dy/dt \bullet \nabla_y \varphi(x, y) \; ds \, dt \, ,$$

*where* $\varphi(\alpha) = \operatorname{csch} \alpha$.

The kernel functions used here have the following significance.

In Gauss's linking integral, the function $\varphi_0(\alpha) = -1/(4\pi\alpha)$, where $\alpha$ is distance from a fixed point, is the fundamental solution of the Laplacian in $R^3$,

$$\Delta\varphi_0 = \delta.$$

Here $\delta$ is the Dirac $\delta$-function.

In formula (1), the function $\varphi_0(\alpha) = (-1/4\pi^2)(\pi - \alpha)\cot\alpha$, is the fundamental solution of the Laplacian on $S^3$,

$$\Delta\varphi_0 = \delta - 1/2\pi^2.$$

Since the volume of $S^3$ is $2\pi^2$, the right-hand side has average value zero.

In formula (2), the function $\varphi(\alpha) = (-1/4\pi^2)(\pi - \alpha)\csc\alpha$ is the fundamental solution of a shifted Laplacian on $S^3$,

$$\Delta\varphi - \varphi = \delta.$$

In formula (3), the function $\varphi(\alpha) = (-1/4\pi)\csch\alpha$ is the fundamental solution of a shifted Laplacian on $H^3$,

$$\Delta\varphi + \varphi = \delta.$$

*Note on typesetting.* In fractions such as $1/4\pi$, $(-1/4\pi^2)$ or $1/2\pi^2$, the $\pi$ and the $\pi^2$ are always in the denominator.

**Comments.**

- For two linked circles in a small ball, the first integral in formula (1) is the dominant term. At the other extreme, for orthogonal great circles, say, the first integral vanishes and the second integral takes the value ±1. In 3-space, the analogue of the second integral would be

$$\int_{K_1 \times K_2} dx/ds \bullet dy/dt \ ds \ dt = \left(\int_{K_1} dx/ds \ ds\right) \bullet \left(\int_{K_2} dy/dt \ dt\right) = 0 \bullet 0 = 0 \ ,$$

and therefore does not appear in the classical Gauss linking integral.

- In formula (2), it makes no difference that parallel transport from $x$ to $y$ on $S^3$ is ambiguous when $y$ is the antipode $-x$ of $x$, because when $\alpha = \pi$ we have $\nabla_y(\varphi(\alpha(x, y))) = 0$.

- The integrand in Gauss's formula, and the integrands in formulas (1), (2) and (3), are invariant under orientation-preserving isometries of the ambient space. In Gauss's formula, this is clear. In formula (1), this follows from the fact that the group of left translations by unit quaternions is a *normal* subgroup of the group SO(4) of all orientation-preserving isometries of $S^3$. In formulas (2) and (3), this is again clear.

- The fact that formula (1) on $S^3$ has no counterpart in $H^3$ is due to the simplicity of the group SO(1, 3) of orientation-preserving isometries of $H^3$.

- Greg Kuperberg, after receiving a copy of the announcement of these results, wrote to us, "As it happens, I thought about this question a few years ago and I obtained one of your formulas, I think, but I never wrote it up." Subsequent correspondence shows that he did indeed obtain, by a beautiful geometric argument totally different from ours, an expression equivalent to formula (2).

## 2. The route to Gauss's linking integral.

According to historian Moritz Epple (1998), Gauss was interested in computing the linking number of the earth's orbit with the orbits of certain asteroids, and although he presented his linking integral formula without proof, it is believed that he simply counted up how many times the vector from the earth to the asteroid covered the "heavenly sphere" ..... a degree-of-map argument. Gauss undoubtedly knew another proof: run a current through the first loop, and calculate the circulation of the resulting magnetic field around the second loop. By Ampere's Law, this circulation is equal to the total current enclosed by the second loop, which means the current flowing along the first loop multiplied by the linking number of the two loops. Then the Biot-Savart formula (1820) for the magnetic field leads directly to Gauss's linking integral.

Gauss's degree-of-map derivation of his linking integral does not work on the 3-sphere $S^3$ because the set of ordered pairs of distinct points in $S^3$ deformation retracts to a 3-sphere rather than to a 2-sphere, as it does in $R^3$. But a degree-of-map derivation, suitably modified, does work in hyperbolic 3-space.

To provide a uniform approach, we develop a steady-state version of classical electrodynamics on $S^3$ and $H^3$, including an explicit formula for the vector-valued Green's operator, an explicit formula of Biot-Savart type for the magnetic field, and a corresponding Ampere's law contained in Maxwell's equations.

Then we follow what we imagine to be Gauss's second line of reasoning, and use the Biot-Savart formula and Ampere's law to derive the linking integrals.

For a development of electrodynamics on bounded subdomains of the 3-sphere, and for applications of the linking, writhing and helicity integrals in this setting, see the Ph.D. thesis of Jason Parsley (2004) . We take pleasure in acknowledging his help throughout the preparation of this paper, and especially in the development of electrodynamics on $S^3$.

We also thank Shea Vick for his critical reading of many parts of this manuscript, and Jozef Dodziuk and Charles Epstein for a number of helpful consultations.

# 3. Magnetic fields in $R^3$, $S^3$ and $H^3$.

In Euclidean 3-space $R^3$, the classical convolution formula of Biot and Savart gives the magnetic field BS(V) generated by a compactly supported current flow V:

$$BS(V)(y) = 1/4\pi \int_{R^3} V(x) \times (y - x) / |y - x|^3 \, dx .$$

For simplicity, we write dx to mean $d(vol_x)$.

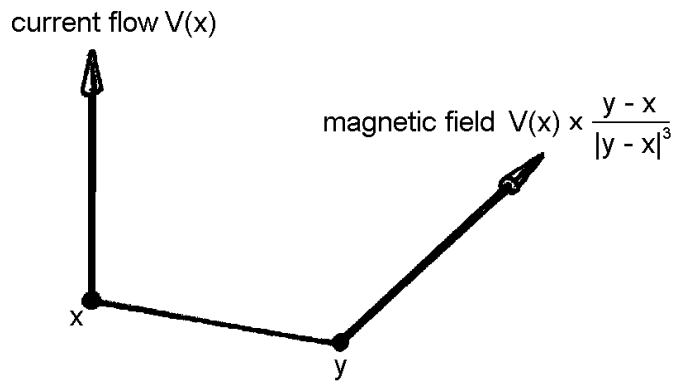

**The current at x contributes to the magnetic field at y**

The Biot-Savart formula can also be written as

$$BS(V)(y) = \int_{R^3} V(x) \times \nabla_y \varphi_0(x, y) \, dx ,$$

where $\varphi_0(\alpha) = -1/(4\pi\alpha)$ is the fundamental solution of the Laplacian in $R^3$.

In $R^3$, if we start with a smooth, compactly supported current flow $V$, then its magnetic field $BS(V)$ is a smooth vector field (although not in general compactly supported) which has the following properties:

(1) It is divergence-free, $\nabla \bullet BS(V) = 0$.

(2) It satisfies Maxwell's equation

$$\nabla_y \times BS(V)(y) = V(y) + \nabla_y \int_{R^3} V(x) \bullet \nabla_x \varphi_0(x, y) \, dx ,$$

where $\varphi_0$ is the fundamental solution of the Laplacian in $R^3$.

(3) $BS(V)(y) \to 0$ as $y \to \infty$.

To see that equation (2) is one of Maxwell's equations, first integrate by parts to get

$$\int_{R^3} V(x) \bullet \nabla_x \varphi_0(x, y) \, dx = - \int_{R^3} (\nabla_x \bullet V(x)) \, \varphi_0(x, y) \, dx ,$$

and then insert this to get

$$\nabla_y \times BS(V)(y) = V(y) - \nabla_y \int_{R^3} (\nabla_x \bullet V(x)) \, \varphi_0(x, y) \, dx .$$

If we think of the vector field $V(x)$ as a steady current, then the negative divergence, $-\nabla_x \bullet V(x)$, is the time rate of accumulation of charge at $x$, and hence the integral

$$-\nabla_y \int_{R^3} (\nabla_x \bullet V(x)) \, \varphi_0(x, y) \, dx$$

is the time rate of increase of the electric field $E$ at $y$. Thus equation (2) is simply Maxwell's equation

$$\nabla \times B = V + \partial E/\partial t .$$

In $R^3$, $S^3$ and $H^3$, a linear operator satisfying conditions (1), (2) and (3) above will be referred to as a ***Biot-Savart operator.***

**Comments.**

- To see that equation (2) above is Maxwell's equation, we integrated by parts, in spite of the fact that the kernel function $\varphi_0(\alpha)$ has a singularity at $\alpha = 0$. We leave it to the reader to check that the validity of this depends on the fact that the singularity is of order $1/\alpha$. We will use this throughout the paper, without further mention.

- Recall Ampere's Law: Given a divergence-free current flow, the circulation of the resulting magnetic field around a loop is equal to the flux of the current through any surface bounded by that loop. This is an immediate consequence of Maxwell's Equation (2) above, since if the current flow $V$ is divergence-free, this equation says that $\nabla \times BS(V) = V$. Then Ampere's Law is just the curl theorem of vector calculus.

In particular, if the current flows along a wire loop, the circulation of the resulting magnetic field around a second loop disjoint from it is equal to the flux of the current through a cross-section of the wire loop, multiplied by the linking number of the two loops. Thus linking numbers are built in to Ampere's Law, and once we have an explicit integral formula for the magnetic field due to a given current flow, we easily get an explicit integral formula for the linking number.

- In $R^3$, conditions (1), (2) and (3) are easily seen to characterize the Biot-Savart operator, as follows. Since conditions (1) and (2) specify the divergence and the curl of $BS(V)$, the difference $BS_1(V) - BS_2(V)$ between two candidates for the Biot-Savart operator would be divergence-free and curl-free. Since $R^3$ is simply connected, this difference would be the gradient of a harmonic function. Hence the components of this gradient must also be harmonic functions. Since they go to zero at infinity, they have to be identically zero. Thus $BS_1(V) = BS_2(V)$.

- In $R^3$, the statement that $BS(V)(y) \to 0$ as $y \to \infty$ can be improved: $BS(V)(y) \to 0$ at infinity like $1/|y|^2$ in general, and if $V$ is divergence-free, then it goes to zero at infinity like $1/|y|^3$. See our paper (2001) on the Biot-Savart operator.

- In $S^3$, conditions (1) and (2) alone suffice to characterize the Biot-Savart operator, since there are no nonzero vector fields on $S^3$ which are simultaneously divergence-free and curl-free.

- In $H^3$, it is not yet clear to us how to characterize the Biot-Savart operator. Even strengthening (3) to require that $BS(V)(y)$ goes to zero exponentially fast at infinity is not quite enough. And in $H^3$, unlike $R^3$, the field $BS(V)$ is not in general of class $L^2$.

# THEOREM 2. BIOT-SAVART INTEGRALS in $S^3$ and $H^3$.

*Biot-Savart operators exist in $S^3$ and $H^3$, and are given by the following formulas, in which V is a smooth, compactly supported vector field:*

**(1)** *On $S^3$ in left-translation format:*

$$BS(V)(y) = \int_{S^3} L_{yx^{-1}} V(x) \times \nabla_y \varphi_0(x, y) \, dx$$
$$- 1/(4\pi^2) \int_{S^3} L_{yx^{-1}} V(x) \, dx$$
$$+ 2 \nabla_y \int_{S^3} L_{yx^{-1}} V(x) \bullet \nabla_y \varphi_1(x, y) \, dx ,$$

*where* $\varphi_0(\alpha) = (-1/4\pi^2)(\pi-\alpha) \cot \alpha$ *and* $\varphi_1(\alpha) = (-1/16\pi^2) \alpha (2\pi-\alpha)$.

**(2)** *On $S^3$ in parallel transport format:*

$$BS(V)(y) = \int_{S^3} P_{yx} V(x) \times \nabla_y \varphi(x, y) \, dx ,$$

*where* $\varphi(\alpha) = (-1/4\pi^2)(\pi - \alpha) \csc \alpha$.

**(3)** *On $H^3$ in parallel transport format:*

$$BS(V)(y) = \int_{H^3} P_{yx} V(x) \times \nabla_y \varphi(x, y) \, dx ,$$

*where* $\varphi(\alpha) = (-1/4\pi) \operatorname{csch} \alpha$.

In formula (1), the function $\varphi_1(\alpha) = (-1/16\pi^2) \alpha (2\pi-\alpha)$ satisfies the equation

$$\Delta \varphi_1 = \varphi_0 - [\varphi_0],$$

where $[\varphi_0]$ denotes the average value of $\varphi_0$ over $S^3$. The other kernel functions already appeared in the linking integrals in Theorem 1.

In formula (3), the magnetic field $BS(V)(y)$ goes to zero at infinity like $e^{-\alpha}$, where $\alpha$ is the distance from y to a fixed point in $H^3$.

# 4. Left-invariant vector fields on $S^3$.

Left-invariant vector fields on $S^3$ play an important role in understanding formula (1) for the Biot-Savart operator. A vector field $V$ on $S^3$ is said to be *left-invariant* if

$$L_{yx^{-1}} V(x) = V(y) \qquad \text{for all } x, y \in S^3.$$

For example, using multiplication of unit quaternions, the vector fields

$$V_1(x) = x\, i, \quad V_2(x) = x\, j, \quad V_3(x) = x\, k$$

form a pointwise orthonormal basis for the 3-dimensional subspace of left-invariant vector fields on $S^3$.

Left-invariant and right-invariant vector fields on $S^3$ are tangent to the orbits of Hopf fibrations of $S^3$ by parallel great circles. Since a steady flow along these parallel great circles is distance-preserving, the left- and right-invariant vector fields are certainly divergence-free.

**PROPOSITION 1.** *Left- and right-invariant vector fields on $S^3$ are curl eigenfields with eigenvalues $-2$ and $+2$, respectively.*

*Proof.* We will show that $\nabla \times V = -2\, V$ for left-invariant fields $V$.

Recall the primitive meaning of "curl". To find the component of $\nabla \times V$ at the point $x$ in a given direction on an oriented Riemannian 3-manifold, you pick a small disk orthogonal at $x$ to that direction, compute the circulation of $V$ around the boundary of this disk, divide by the area of the disk, and then take the limit as the disk shrinks to the point $x$. To keep signs straight, the orientation of the disk, combined with the orientation of the normal direction, should give the chosen orientation on the ambient space. In this view, the primitive meaning of "curl" is given by the "curl theorem", i.e., Stokes' theorem.

To apply this to left-invariant vector fields on $S^3$, first note that any two of these are equivalent to one another by an orientation-preserving rigid motion of $S^3$, so it is enough to prove the result for one such field.

For example, if

$$x = x_0 + x_1 i + x_2 j + x_3 k,$$

we will use

$$V(x) = x i = (x_0 + x_1 i + x_2 j + x_3 k) i$$

$$= -x_1 + x_0 i + x_3 j - x_2 k.$$

Furthermore, $V$ itself is invariant under a transitive group of orientation-preserving rigid motions of $S^3$, so it is enough to show that $\nabla \times V = -2\,V$ at a single point. Choose the point $x = 1$, where $V(x) = i$. Orient $S^3$ as usual, so that $i, j, k$ is an ordered basis for the tangent space at the point $x = 1$. Now we find the component of $\nabla \times V$ at the point $x = 1$ in the i-direction, as follows.

To get a small disk through $x$ orthogonal to this direction, use the spherical cap of geodesic radius $\alpha$ on the great 2-sphere in the 1jk-subspace. The typical point along the boundary of this disk is

$$y = \cos \alpha + (\sin \alpha \cos \varphi) j + (\sin \alpha \sin \varphi) k.$$

Note that

$$dy/d\varphi = (-\sin \alpha \sin \varphi) j + (\sin \alpha \cos \varphi) k,$$

and

$$V(y) = y i = (\cos \alpha) i + (\sin \alpha \sin \varphi) j - (\sin \alpha \cos \varphi) k.$$

The circulation of $V$ around the boundary of this disk is given by

$$\int_0^{2\pi} V(y) \bullet dy/d\varphi \; d\varphi \;=\; \int_0^{2\pi} -\sin^2\alpha \; d\varphi \;=\; -2\pi \sin^2\alpha.$$

The area of this disk is given by

$$\int_0^\alpha 2\pi \sin \alpha \; d\alpha \;=\; 2\pi (1 - \cos \alpha).$$

Dividing circulation by area, we get

$$-\sin^2\alpha / (1 - \cos \alpha) \;=\; -(1 - \cos^2\alpha)/(1 - \cos \alpha) \;=\; -(1 + \cos \alpha).$$

As we let $\alpha \to 0$, this approaches $-2$. Hence at the point $x = 1$, the component of $\nabla \times V$ in the i-direction is $-2$. So far, so good, since $V(x) = i$.

A similar computation shows that the components of $\nabla \times V$ in the j- and k-directions are zero. Thus, at the point $x = 1$, we have $\nabla \times V = -2\,i = -2\,V$. As mentioned earlier, this implies that $\nabla \times V = -2\,V$ at every point of $S^3$, completing our proof.

Likewise, if V is right-invariant, then $\nabla \times V = +2\,V$.

We noted in section 3 that if V is divergence-free, then $\nabla \times BS(V) = V$. If V is left-invariant, then we also have $\nabla \times (-\tfrac{1}{2}V) = V$. Therefore $BS(V)$ and $-\tfrac{1}{2}V$ are divergence-free fields with the same curl. Their difference is a field which is simultaneously divergence-free and curl-free, and hence, since $S^3$ is simply connected, must be zero. Thus $BS(V) = -\tfrac{1}{2}V$.

Likewise, if V is right-invariant, then $BS(V) = +\tfrac{1}{2}V$.

Let $VF(S^3)$ denote the set of all smooth vector fields on $S^3$, regarded as an infinite-dimensional vector space with the $L^2$ inner product

$$<V, W> = \int_{S^3} V \bullet W \ d(\text{vol}).$$

The following comments, which we do not prove here, help to explain the different roles played by the three integrals in formula (1) of Theorem 2 above.

(a) *The first integral is zero if and only if V is a left-invariant field.*

(b) *The second integral is zero if and only if V is orthogonal, in the above inner product, to the 3-dimensional subspace of left-invariant fields.*

(c) *The third integral is zero if and only if V is divergence-free.*

(d) *If V is a nonzero gradient field, then the first and third integrals are nonzero and cancel each other, while the second integral is zero.*

# 5. Scalar-valued Green's operators in $R^3$, $S^3$ and $H^3$.

The scalar Laplacian on any Riemannian manifold is given by $\Delta f = \nabla \cdot \nabla f$.

Suppose that f is a smooth real-valued function on $R^3$, $S^3$ or $H^3$ which depends only on the distance $\alpha$ from some fixed point. Then

$$\Delta f = \frac{1}{\alpha^2} \frac{d}{d\alpha}\left(\alpha^2 \frac{df}{d\alpha}\right) \qquad \text{in } R^3,$$

$$\Delta f = \frac{1}{\sin^2\alpha} \frac{d}{d\alpha}\left(\sin^2\alpha \frac{df}{d\alpha}\right) \qquad \text{on } S^3,$$

$$\Delta f = \frac{1}{\sinh^2\alpha} \frac{d}{d\alpha}\left(\sinh^2\alpha \frac{df}{d\alpha}\right) \qquad \text{in } H^3.$$

Let f be a smooth real-valued function with compact support on $R^3$, $S^3$ or $H^3$. Then the *scalar-valued Green's operators* on these three spaces are given by the convolution formulas

$$Gr(f)(y) = \int f(x)\, \varphi_0(x, y)\, dx,$$

where the integral is taken over the entire space, and where $\varphi_0$ is the fundamental solution of the Laplacian, that is,

$$\varphi_0(\alpha) = -1/(4\pi\alpha) \qquad \text{in } R^3,$$

$$\varphi_0(\alpha) = (-1/4\pi^2)\,(\pi - \alpha)\cot\alpha \qquad \text{on } S^3,$$

$$\varphi_0(\alpha) = (-1/4\pi)\,(\coth\alpha - 1) \qquad \text{in } H^3.$$

In $R^3$ and $H^3$ we have

(1) $\Delta \, Gr(f) = f$, and

(2) $Gr(f)(y) \to 0$ as $y \to \infty$.

On $S^3$ we have

(1') $\Delta \, Gr(f) = f - [f]$, and

(2') $[Gr(f)] = 0$,

where $[f]$ denotes the average value of $f$ over $S^3$.

**Comments.**

- In $R^3$, conditions (1) and (2) characterize the Green's operator, because the difference $Gr_1(f) - Gr_2(f)$ between two candidates for the Green's operator is a harmonic function which goes to zero at infinity, and hence must be identically zero.

- In $S^3$, conditions (1') and (2') characterize the Green's operator, because the scalar Laplacian is a bijection from the set of smooth functions with average value zero to itself.

- In $H^3$, it is not clear to us how to characterize the scalar-valued Green's operator. Even strengthening (2) to require that $Gr(f)(y)$ goes to zero exponentially fast at infinity is not quite enough.

# 6. Vector-valued Green's operators in $R^3$, $S^3$ and $H^3$.

The *vector Laplacian* in $R^3$, $S^3$ and $H^3$,

$$\Delta V = -\nabla \times \nabla \times V + \nabla(\nabla \bullet V),$$

corresponds to the Laplacian on 1-forms,

$$\Delta \beta = (d\delta + \delta d)\beta,$$

where $\delta = \pm * d *$, with $*$ the Hodge star operator, in the following sense. Given the vector field $V$, if we define the 1-form $\downarrow V$ by $(\downarrow V)(W) = V \bullet W$, then we get $\Delta(\downarrow V) = \downarrow(\Delta V)$. The symbol $\downarrow$ reminds us of lowering indices in tensor notation.

In $R^3$, the vector Laplacian can be computed component-wise in rectangular coordinates.

If $V$ is a smooth compactly supported vector field in $R^3$, then the *Green's operator* acting on $V$ is given by the convolution formula

$$Gr(V)(y) = \int_{R^3} V(x) \, \varphi_0(x, y) \, dx,$$

where $\varphi_0(\alpha) = -1/(4\pi\alpha)$ is the fundamental solution of the Laplacian. It has the following properties:

(1) $\Delta \, Gr(V) = V$.

(2) $Gr(V)(y) \to 0$ as $y \to \infty$.

In $R^3$, $S^3$ and $H^3$, a linear operator satisfying conditions (1) and (2) above will be referred to as a *vector-valued Green's operator.*

**Comments.**

- In $R^3$, conditions (1) and (2) are easily seen to characterize the Green's operator, as follows. The difference $Gr_1(V) - Gr_2(V)$ between two candidates for the Green's operator is a harmonic vector field, and therefore its three components in rectangular coordinates are harmonic functions. By property (2), they go to zero at infinity and hence must be identically zero. Thus $Gr_1(V) = Gr_2(V)$.

- In $R^3$, the statement that $Gr(V)(y) \to 0$ as $y \to \infty$ can be improved: $Gr(V)(y) \to 0$ at infinity like $1/|y|$ in general, and if $V$ is divergence-free, then it goes to zero at infinity like $1/|y|^2$. See our paper (2001) on the Biot-Savart operator.

- In $S^3$, condition (1) alone suffices to characterize the Green's operator, since the vector Laplacian $\Delta: VF(S^3) \to VF(S^3)$ is bijective.

- In $S^3$, because the vector Laplacian is bijective and commutes with grad, curl and div, we immediately get the same for the scalar and vector-valued Green's operators:

$$\nabla \, Gr(f) = Gr(\nabla f) \, , \quad \nabla \times Gr(V) = Gr(\nabla \times V) \, , \quad \nabla \bullet Gr(V) = Gr(\nabla \bullet V) \, .$$

- In $H^3$, it is not yet clear to us how to characterize the vector-valued Green's operator. Even strengthening (2) to require that $Gr(V)(y)$ goes to zero exponentially fast at infinity is not quite enough.

The following theorem asserts the existence of vector-valued Green's operators in $S^3$ and $H^3$, and provides explicit integral formulas.

# THEOREM 3. GREEN'S OPERATORS in $S^3$ and $H^3$.

*Vector-valued Green's operators exist in $S^3$ and $H^3$, and are given by the following formulas, in which $V$ is a smooth, compactly supported vector field:*

**(1)** *On $S^3$ in left-translation format:*

$$Gr(V)(y) = \int_{S^3} L_{yx^{-1}} V(x)\, \varphi_0(x, y)\, dx$$

$$+ \; 2 \int_{S^3} L_{yx^{-1}} V(x) \times \nabla_y \varphi_1(x, y)\, dx$$

$$+ \; 4 \nabla_y \int_{S^3} L_{yx^{-1}} V(x) \bullet \nabla_y \varphi_2(x, y)\, dx\,,$$

*where*

$$\varphi_0(\alpha) = (-1/4\pi^2)\, (\pi - \alpha) \cot \alpha$$

$$\varphi_1(\alpha) = (-1/16\pi^2)\, \alpha\, (2\pi - \alpha)$$

$$\varphi_2(\alpha) = (-1/192\pi^2)\, \bigl(3\alpha(2\pi - \alpha) + 2\alpha(\pi - \alpha)(2\pi - \alpha)\cot \alpha\bigr),$$

*and these three kernel functions are related by*

$$\varphi_2 \xrightarrow{\Delta} \varphi_1 - [\varphi_1] \xrightarrow{\Delta} \varphi_0 - [\varphi_0] \xrightarrow{\Delta} \delta - [\delta]\,.$$

**(2)** *On $S^3$ in parallel transport format:*

$$\mathbf{Gr(V)(y) \;=\; \int_{S^3} P_{yx}V(x)\; \varphi_2(x, y)\; dx \;+\; \nabla_y\!\int_{S^3} P_{yx}V(x)\bullet \nabla_y\varphi_3(x, y)\; dx,}$$

*where*

$$\varphi_2(\alpha) \;=\; (-1/4\pi^2)\,(\pi - \alpha)\csc\alpha \;+\; (1/8\pi^2)\,(\pi - \alpha)^2 / (1 + \cos\alpha)$$

*and*

$$\varphi_3(\alpha) \;=\; (-1/24)\,(\pi - \alpha)\cot\alpha \;-\; (1/16\pi^2)\,\alpha\,(2\pi - \alpha)$$

$$+\; (1/8\pi^2) \int_\alpha^\pi \left( (\pi - \alpha)^3 / (3\sin^2\alpha) \right) + \left( (\pi - \alpha)^2 / \sin\alpha \right) \, d\alpha.$$

**(3)** *On $H^3$ in parallel transport format:*

$$\mathbf{Gr(V)(y) \;=\; \int_{S^3} P_{yx}V(x)\; \varphi_2(x, y)\; dx \;+\; \nabla_y\!\int_{S^3} P_{yx}V(x)\bullet \nabla_y\varphi_3(x, y)\; dx,}$$

*where*

$$\varphi_2(\alpha) \;=\; (-1/4\pi)\,\operatorname{csch}\alpha \;+\; (1/4\pi)\,\alpha / (1 + \cosh\alpha)$$

*and*

$$\varphi_3(\alpha) \;=\; (1/4\pi)\,\alpha / (e^{2\alpha} - 1) \;+\; (1/4\pi)\int_0^\alpha \left( (\alpha / \sinh\alpha) - (\alpha^2 / 2\sinh^2\alpha) \right) d\alpha.$$

**Comment.** In formula (2), the third term in the expression for the kernel function $\varphi_3(\alpha)$ is left as a definite integral because it is a non-elementary function involving polylogarithms of indices 2 and 3, where the polylogarithm of index s is defined by

$$\text{polylog}(s, x) \;=\; \text{Li}_s(x) \;=\; \sum_{n\geq 1} x^n / n^s.$$

In formula (3), the second term in the expression for the kernel function $\varphi_3(\alpha)$ is left as a definite integral because it is a non-elementary function involving dilogarithms, that is, polylogarithms of index 2.

# 7. Classical electrodynamics and Maxwell's equations on $S^3$ and $H^3$.

As mentioned in section 2, our route to the Gauss linking integral begins by developing a steady-state version of classical electrodynamics.

Maxwell's equations, so familiar to us in Euclidean 3-space $R^3$, can also be shown to hold on $S^3$ and $H^3$, with the appropriate definitions, as follows.

Let $\rho$ be a smooth real-valued function, which we think of as a charge density. On $R^3$ and $H^3$, we ask that $\rho$ have compact support, while on $S^3$ we ask that it have average value zero.

In each case, we define the corresponding electric field $E = E(\rho)$ by the formula

$$E(\rho)(y) = \nabla_y \int \rho(x) \varphi_0(x, y) \, dx \,,$$

where $\varphi_0$ is the fundamental solution of the Laplacian.

Let $V$ be a smooth vector field, which we think of as a steady-state (i.e., time-independent) current distribution. On $R^3$ and $H^3$, we ask that $V$ have compact support, while on $S^3$ we impose no restriction on $V$.

In each case, we define the corresponding magnetic field $B = BS(V)$ by one of the formulas in Theorem 2. Since $V$ is steady-state, so is $B$, and therefore $\partial B/\partial t = 0$.

**THEOREM 4. MAXWELL'S EQUATIONS in $S^3$ and $H^3$.**

*With these definitions, Maxwell's equations,*

$$\nabla \cdot E = \rho \qquad \nabla \times E = 0$$

$$\nabla \cdot B = 0 \qquad \nabla \times B = V + \partial E/\partial t \,,$$

*hold on $S^3$ and in $H^3$, just as they do in $R^3$.*

Additionally, as in $R^3$, Ampere's Law follows immediately from the last of the Maxwell equations and the curl theorem.

# B.  PROOFS ON $S^3$ IN LEFT-TRANSLATION FORMAT

## 8.  Proof scheme for Theorem 3, formula (1).

When working on $S^3$ in left-translation format, our first order of business will be to derive formula (1) of Theorem 3 for the vector-valued Green's operator.

It is natural to hope, by analogy with $R^3$, that

$$A(V, \varphi_0)(y) \;=\; \int_{S^3} L_{yx^{-1}} V(x)\, \varphi_0(x, y)\, dx$$

will be the Green's operator on $S^3$ when $\varphi_0$ is the fundamental solution of the scalar Laplacian.

It doesn't turn out this way, but the fundamental solution of the Laplacian is not the only possible kernel function, and $A(V, \varphi)$ is not the only way of convolving a vector field with a kernel to get another vector field.

Consider the convolutions

$$A(V, \varphi)(y) \;=\; \int_{S^3} L_{yx^{-1}} V(x)\, \varphi(x, y)\, dx$$

$$B(V, \varphi)(y) \;=\; \int_{S^3} L_{yx^{-1}} V(x) \times \nabla_y\, \varphi(x, y)\, dx$$

$$G(V, \varphi)(y) \;=\; \nabla_y \int_{S^3} L_{yx^{-1}} V(x) \bullet \nabla_y\, \varphi(x, y)\, dx \,.$$

We will develop a "calculus of vector convolutions", in which we compute formulas for curls, divergences and Laplacians of convolutions, first in $R^3$ to serve as a model, and then on $S^3$ in left-translation format, and in particular obtain

**PROPOSITION 2.** *On $S^3$, the vector Laplacian has the following effect on convolutions of types* **A**, **B** *and* **G**:

$$\Delta A(V, \varphi) = A(V, \Delta\varphi) \quad - 4\, A(V, \varphi) \quad - 2\, B(V, \varphi)$$

$$\Delta B(V, \varphi) = B(V, \Delta\varphi) \quad + 2\, A(V, \Delta\varphi) - 2\, G(V, \varphi)$$

$$\Delta G(V, \varphi) = G(V, \Delta\varphi) \,.$$

The large spaces after the first terms on the right hand sides of the A and B lines serve as a reminder that in Euclidean space $R^3$, only those first terms appear. In particular, in $R^3$ we have $\Delta A(V, \varphi) = A(V, \Delta\varphi)$. So naturally, when we choose $\varphi = \varphi_0$ to be a fundamental solution of the scalar Laplacian in $R^3$, we get $\Delta A(V, \varphi_0) = A(V, \Delta\varphi_0) = A(V, \delta) = V$, which tells us that $A(V, \varphi_0) = Gr(V)$ there.

Once we have Proposition 2 in hand, we will then get formula (1) of Theorem 3 for the Green's operator on $S^3$ by adding one copy of line A to two copies of line B and four copies of line G, each with a different choice of kernel function $\varphi$.

# 9. The calculus of vector convolutions in $R^3$.

Let $\alpha = |x - y|$ denote the distance between the points $x$ and $y$ in $R^3$, and let $\varphi(\alpha)$ be some unspecified function of $\alpha$. As usual, we'll write $\varphi(x, y)$ in place of $\varphi(\alpha(x, y))$ for simplicity.

Let $V$ be a smooth vector field defined on $R^3$ and having compact support. Keeping the function $\varphi$ unspecified, we define vector fields $A(V, \varphi)$ and $B(V, \varphi)$, a scalar function $g(V, \varphi)$ and another vector field $G(V, \varphi)$ as follows.

$$A(V, \varphi)(y) = \int_{R^3} V(x)\, \varphi(x, y)\, dx$$

$$B(V, \varphi)(y) = \int_{R^3} V(x) \times \nabla_y \varphi(x, y)\, dx$$

$$g(V, \varphi)(y) = \int_{R^3} V(x) \bullet \nabla_y \varphi(x, y)\, dx$$

$$G(V, \varphi)(y) = \nabla_y g(V, \varphi)(y) = \nabla_y \int_{R^3} V(x) \bullet \nabla_y \varphi(x, y)\, dx\ .$$

If we use the fundamental solution $\varphi_0(\alpha) = -1/(4\pi\alpha)$ of the scalar Laplacian as the kernel in the above convolutions, then $A(V, \varphi_0)$ gives the Green's operator, and if we think of $V$ as a current flow, then $B(V, \varphi_0)$ is the resulting magnetic field, $-A(V, \varphi_0)$ is its vector potential, and $-G(V, \varphi_0)$ is the rate of change $\partial E/\partial t$ of the electric field, caused by the accumulation and dissipation of charge by the current flow $V$.

The curl and divergence of the vector fields $A(V, \varphi)$ and $B(V, \varphi)$ are given by

**LEMMA 1.** $\quad \nabla \times A(V, \varphi) = - B(V, \varphi) \qquad\qquad \nabla \bullet A(V, \varphi) = g(V, \varphi)$

$$\nabla \times B(V, \varphi) = A(V, \Delta\varphi) - G(V, \varphi) \qquad \nabla \bullet B(V, \varphi) = 0 \, .$$

**Proof.** The proofs are exercises in differentiating under the integral sign and then using the various Leibniz rules from vector calculus. We will show the brief arguments for the first and last of the four formulas above, since these two will change when we move to the 3-sphere. The other two formulas, whose proofs we leave to the reader, are the same in $R^3$ and $S^3$.

To begin with the first formula,

$$\nabla_y \times A(V, \varphi)(y) = \int_{R^3} \nabla_y \times \{V(x) \, \varphi(x, y)\} \, dx$$

$$= \int_{R^3} (\nabla_y \times V(x)) \, \varphi(x, y) - V(x) \times \nabla_y \varphi(x, y) \, dx$$

$$= - \int_{R^3} V(x) \times \nabla_y \varphi(x, y) \, dx$$

$$= - B(V, \varphi)(y) \, ,$$

as claimed, since $\nabla_y \times V(x) = 0$. When we move to $S^3$, the corresponding term will *not* be zero.

As for the last formula,

$$\nabla_y \bullet B(V, \varphi)(y) = \int_{R^3} \nabla_y \bullet \{V(x) \times \nabla_y \varphi(x, y)\} \, dx$$

$$= \int_{R^3} (\nabla_y \times V(x)) \bullet \nabla_y \varphi(x, y) - V(x) \bullet (\nabla_y \times \nabla_y \varphi(x, y)) \, dx$$

$$= 0 \, ,$$

as claimed, because $\nabla_y \times V(x) = 0$ and $\nabla_y \times \nabla_y \varphi(x, y) = 0$. Again, when we move to $S^3$, the expression corresponding to $\nabla_y \times V(x)$ will not be zero.

**Comment.**

If we use the fundamental solution $\varphi_0$ of the Laplacian as our kernel, think of V as a current flow and $B(V, \varphi_0)$ as its magnetic field, we get

$$\nabla \times B(V, \varphi_0) = A(V, \Delta\varphi_0) - G(V, \varphi_0),$$

which is Maxwell's equation,

$$\nabla \times B = V + \partial E/\partial t.$$

Now we turn to the Laplacians of our convolutions.

**LEMMA 2.** $\qquad \Delta A(V, \varphi) = A(V, \Delta\varphi) \qquad\qquad \Delta B(V, \varphi) = B(V, \Delta\varphi)$

$\qquad\qquad\qquad \Delta g(V, \varphi) = g(V, \Delta\varphi) \qquad\qquad \Delta G(V, \varphi) = G(V, \Delta\varphi).$

**Proof.**

The arguments are straightforward, using the formula for the vector Laplacian, the information from Lemma 1, and the fact that Laplacian and gradient commute. We write it out for the first formula, and leave the other three to the reader.

$$\Delta A(V, \varphi) = -\nabla \times \nabla \times A(V, \varphi) + \nabla(\nabla \bullet A(V, \varphi))$$

$$= -\nabla \times (-B(V, \varphi)) + \nabla g(V, \varphi)$$

$$= A(V, \Delta\varphi) - G(V, \varphi) + G(V, \varphi)$$

$$= A(V, \Delta\varphi).$$

# 10. The calculus of vector convolutions in $S^3$.

We turn now to the calculus of vector convolutions in $S^3$. Our specific goal is to prove Theorem 3, formula (1).

Let $V$ be a smooth vector field defined on $S^3$. For any kernel function $\varphi(\alpha)$, we define

$$A(V, \varphi)(y) = \int_{S^3} L_{yx^{-1}} V(x) \, \varphi(x, y) \, dx$$

$$B(V, \varphi)(y) = \int_{S^3} L_{yx^{-1}} V(x) \times \nabla_y \varphi(x, y) \, dx$$

$$g(V, \varphi)(y) = \int_{S^3} L_{yx^{-1}} V(x) \bullet \nabla_y \varphi(x, y) \, dx$$

$$G(V, \varphi)(y) = \nabla_y \, g(V, \varphi)(y) = \nabla_y \int_{S^3} L_{yx^{-1}} V(x) \bullet \nabla_y \varphi(x, y) \, dx \, .$$

We now continue as in $R^3$: the curl and divergence of the vector fields $A(V, \varphi)$ and $B(V, \varphi)$ are given by

**LEMMA 3.**

$$\nabla \times A(V, \varphi) = -B(V, \varphi) \qquad\qquad - 2 A(V, \varphi)$$

$$\nabla \bullet A(V, \varphi) = g(V, \varphi)$$

$$\nabla \times B(V, \varphi) = A(V, \Delta\varphi) - G(V, \varphi)$$

$$\nabla \bullet B(V, \varphi) = \qquad\qquad\qquad\qquad - 2 g(V, \varphi) \, .$$

We have displayed these formulas so that the terms which already appeared in the corresponding formulas in $R^3$, given in Lemma 1, are to the left, while those which are new in $S^3$ are to the far right.

**Proof.** We give the arguments for the first and last of these formulas, each of which looks like the corresponding formula in $R^3$ with an extra term added, and leave the other two formulas to the reader.

We begin with the first formula,

$$\nabla_y \times A(V, \varphi)(y) = \nabla_y \times \int_{S^3} L_{yx^{-1}} V(x) \, \varphi(x, y) \, dx$$

$$= \int_{S^3} \nabla_y \times \{L_{yx^{-1}} V(x) \, \varphi(x, y)\} \, dx$$

$$= \int_{S^3} \left(\nabla_y \times L_{yx^{-1}} V(x)\right) \varphi(x, y) - L_{yx^{-1}} V(x) \times \nabla_y \, \varphi(x, y) \, dx \, .$$

In $R^3$, $\nabla_y \times V(x)$ was zero in the first term because we were differentiating with respect to $y$ a vector field $V(x)$ which did not depend on $y$.

In $S^3$, the vector field $L_{yx^{-1}} V(x)$ *does* depend on $y$, and its curl, $\nabla_y \times L_{yx^{-1}} V(x)$, is *not* zero. In fact, as we indicated in section 4, the vector field $L_{yx^{-1}} V(x)$ is left-invariant with respect to $y$, and as such is a curl-eigenfield with eigenvalue $-2$, that is,

$$\nabla_y \times L_{yx^{-1}} V(x) = -2 L_{yx^{-1}} V(x).$$

This observation lets us continue and complete the calculation:

$$\nabla_y \times A(V, \varphi)(y) = \int_{S^3} [\{\nabla_y \times L_{yx^{-1}} V(x)\} \varphi(x, y) - L_{yx^{-1}} V(x) \times \nabla_y \varphi(x, y)] \, dx$$

$$= \int_{S^3} [-2 L_{yx^{-1}} V(x) \varphi(x, y) - L_{yx^{-1}} V(x) \times \nabla_y \varphi(x, y)] \, dx$$

$$= -2 \int_{S^3} L_{yx^{-1}} V(x) \varphi(x, y) \, dx - \int_{S^3} L_{yx^{-1}} V(x) \times \nabla_y \varphi(x, y) \, dx$$

$$= -2 A(V, \varphi)(y) - B(V, \varphi)(y),$$

as claimed.

As for the last formula,

$$\nabla_y \bullet B(V, \varphi)(y) = \nabla_y \bullet \int_{S^3} L_{yx^{-1}} V(x) \times \nabla_y \varphi(x, y) \, dx$$

$$= \int_{S^3} \nabla_y \bullet \{L_{yx^{-1}} V(x) \times \nabla_y \varphi(x, y)\} \, dx$$

$$= \int_{S^3} [\{\nabla_y \times L_{yx^{-1}} V(x)\} \bullet \nabla_y \varphi(x, y) \, dx - L_{yx^{-1}} V(x) \bullet \{\nabla_y \times \nabla_y \varphi(x, y)\}] \, dx$$

$$= \int_{S^3} -2 L_{yx^{-1}} V(x) \bullet \nabla_y \varphi(x, y) \, dx$$

$$= -2 g(V, \varphi)(y),$$

as claimed, since $\nabla_y \times \nabla_y \varphi(x, y) = 0$.

**Proof of Proposition 2.**

$$\Delta A(V, \varphi) = -\nabla \times \nabla \times A(V, \varphi) + \nabla(\nabla \bullet A(V, \varphi))$$

$$= -\nabla \times (-B(V, \varphi) - 2 A(V, \varphi)) + \nabla g(V, \varphi)$$

$$= \nabla \times B(V, \varphi) + 2 \nabla \times A(V, \varphi) + G(V, \varphi)$$

$$= A(V, \Delta\varphi) - G(V, \varphi) - 2 B(V, \varphi) - 4 A(V, \varphi) + G(V, \varphi)$$

$$= A(V, \Delta\varphi) - 4 A(V, \varphi) - 2 B(V, \varphi),$$

as claimed.

$$\Delta B(V, \varphi) = -\nabla \times \nabla \times B(V, \varphi) + \nabla(\nabla \bullet B(V, \varphi))$$

$$= -\nabla \times (A(V, \Delta\varphi) - G(V, \varphi)) + \nabla(-2 g(V, \varphi))$$

$$= -\nabla \times A(V, \Delta\varphi) - 2 G(V, \varphi)$$

$$= B(V, \Delta\varphi) + 2 A(V, \Delta\varphi) - 2 G(V, \varphi),$$

as claimed.

$$\Delta g(V, \varphi) = \nabla \bullet \nabla g(V, \varphi)$$

$$= \nabla \bullet G(V, \varphi)$$

$$= \nabla \bullet \{A(V, \Delta\varphi) - \nabla \times B(V, \varphi)\}$$

$$= g(V, \Delta\varphi),$$

using the formulas for $\nabla \times B(V, \varphi)$ and for $\nabla \bullet A(V, \Delta\varphi)$ obtained in Lemma 3.

Finally,

$$\Delta G(V, \varphi) = \Delta \nabla g(V, \varphi)$$

$$= \nabla \Delta g(V, \varphi)$$

$$= \nabla g(V, \Delta\varphi)$$

$$= G(V, \Delta\varphi),$$

completing the proof of Proposition 2.

## 11. Proof of Theorem 3, formula (1).

We intend to show that on $S^3$ in left-translation format, the Green's operator on vector fields is given by

$$\text{Gr}(V)(y) = \int_{S^3} L_{yx^{-1}} V(x) \, \varphi_0(x, y) \, dx$$

$$+ \, 2 \int_{S^3} L_{yx^{-1}} V(x) \times \nabla_y \varphi_1(x, y) \, dx$$

$$+ \, 4 \nabla_y \int_{S^3} L_{yx^{-1}} V(x) \bullet \nabla_y \varphi_2(x, y) \, dx$$

$$= A(V, \varphi_0)(y) + 2 B(V, \varphi_1)(y) + 4 G(V, \varphi_2)(y) \, ,$$

where

$$\varphi_0(\alpha) = (-1/4\pi^2)(\pi - \alpha) \cot \alpha$$

$$\varphi_1(\alpha) = (-1/16\pi^2) \, \alpha \, (2\pi - \alpha)$$

$$\varphi_2(\alpha) = (-1/192\pi^2) \left( 3\alpha(2\pi - \alpha) + 2\alpha(\pi - \alpha)(2\pi - \alpha) \cot \alpha \right) \, .$$

Direct computation shows that these three kernel functions are related by

$$\varphi_2 \xrightarrow{\Delta} \varphi_1 - [\varphi_1] \xrightarrow{\Delta} \varphi_0 - [\varphi_0] \xrightarrow{\Delta} \delta - [\delta] \, ,$$

with average values

$$[\varphi_0] = -1/(8\pi^2) \quad \text{and} \quad [\varphi_1] = -1/(32\pi^2) - 1/24 \, .$$

To show that $\Delta Gr(V) = V$, we add up the three terms in $\Delta Gr(V)$ as follows.

$$\Delta A(V, \varphi_0) = A(V, \Delta\varphi_0) - 4\,A(V, \varphi_0) - 2\,B(V, \varphi_0)$$
$$= A(V, \Delta\varphi_0) - 4\,A(V, \varphi_0) - 2\,B(V, \varphi_0 - [\varphi_0])$$
$$2\,\Delta B(V, \varphi_1) = \qquad\qquad 4\,A(V, \Delta\varphi_1) + 2\,B(V, \Delta\varphi_1) - 4\,G(V, \varphi_1)$$
$$= \qquad\qquad 4\,A(V, \Delta\varphi_1) + 2\,B(V, \Delta\varphi_1) - 4\,G(V, \varphi_1 - [\varphi_1])$$
$$4\,\Delta G(V, \varphi_2) = \qquad\qquad\qquad\qquad 4\,G(V, \Delta\varphi_2)\,.$$

Notice that we take the gradient of the kernel function $\varphi_0$ in the definition of $B(V, \varphi_0)$, and so replacing $\varphi_0$ by $\varphi_0 - [\varphi_0]$ causes no change in value. Likewise, we replace $\varphi_1$ by $\varphi_1 - [\varphi_1]$ in the definition of $G(V, \varphi_1)$ without changing its value.

Now we add up the columns, use the facts that $\Delta\varphi_1 = \varphi_0 - [\varphi_0]$ and $\Delta\varphi_2 = \varphi_1 - [\varphi_1]$, and get

$$\Delta Gr(V) = \Delta A(V, \varphi_0) + 2\,\Delta B(V, \varphi_1) + 4\,\Delta G(V, \varphi_2)$$
$$= A(V, \Delta\varphi_0) - 4\,A(V, \varphi_0) + 4\,A(V, \Delta\varphi_1)$$
$$= A(V, \delta - [\delta]) - 4\,A(V, \varphi_0) + 4\,A(V, \varphi_0 - [\varphi_0])$$
$$= A(V, \delta) - A(V, [\delta]) - 4\,A(V, [\varphi_0])\,.$$

Note that $[\delta] = 1/2\pi^2$, and $[\varphi_0] = -1/8\pi^2$, so the last two terms above cancel, and we get

$$\Delta Gr(V) = A(V, \delta) = V\,,$$

completing the proof of Theorem 3, formula (1).

## 12. Proof of Theorem 2, formula (1).

Let $V$ be a smooth vector field on $S^3$. Thinking of $V$ as a steady current flow, we define the corresponding magnetic field $BS(V)$ by the formula

$$BS(V) = -\nabla \times Gr(V).$$

It is easy to check that $BS(V)$ satisfies the two properties required of it in section 4, as follows.

(1) $\nabla \bullet BS(V) = 0$ because the divergence of a curl is always zero.

(2) To get Maxwell's equation, we compute:

$$\nabla \times BS(V)(y) = -\nabla \times \nabla \times Gr(V)(y)$$

$$= \Delta Gr(V)(y) - \nabla(\nabla \bullet Gr(V)(y)$$

$$= V(y) - \nabla(Gr(\nabla \bullet V)(y)$$

$$= V(y) - \nabla_y \int_{S^3} (\nabla_x \bullet V(x)) \, \varphi_0(x, y) \, dx$$

$$= V(y) + \nabla_y \int_{S^3} V(x) \bullet \nabla_x \varphi_0(x, y) \, dx \, .$$

where $\varphi_0$ is the fundamental solution of the scalar Laplacian on $S^3$. Note that we used the fact, mentioned in section 6, that the Green's operator on $S^3$ commutes with divergence.

The explicit expression for BS(V) given in formula (1) of Theorem 2 comes from the explicit formula for Gr(V) in left-translation format, given in Theorem 3, formula (1), together with the curl formulas given in Lemma 3, as follows.

$$BS(V) = -\nabla \times Gr(V)$$

$$= -\nabla \times A(V, \varphi_0) - 2\nabla \times B(V, \varphi_1) - 4\nabla \times G(V, \varphi_2)$$

$$= B(V, \varphi_0) + 2A(V, \varphi_0) - 2A(V, \Delta\varphi_1) + 2G(V, \varphi_1)$$

$$= B(V, \varphi_0) + 2A(V, \varphi_0) - 2A(V, \varphi_0 - [\varphi_0]) + 2G(V, \varphi_1)$$

$$= B(V, \varphi_0) + 2A(V, [\varphi_0]) + 2G(V, \varphi_1)$$

$$= B(V, \varphi_0) - (1/4\pi^2)A(V, 1) + 2G(V, \varphi_1)$$

$$= \int_{S^3} L_{yx^{-1}} V(x) \times \nabla_y \varphi_0(x, y) \, dx$$

$$- 1/(4\pi^2) \int_{S^3} L_{yx^{-1}} V(x) \, dx$$

$$+ 2 \nabla_y \int_{S^3} L_{yx^{-1}} V(x) \bullet \nabla_y \varphi_1(x, y) \, dx ,$$

as claimed.

# C. PROOFS ON $S^3$ IN PARALLEL TRANSPORT FORMAT

## 13. Parallel transport in $S^3$.

We want to switch formats now, from left-translation to parallel transport. First we explain how to conduct vector calculus in this mode. Then we derive formula (2) of Theorem 2 for the Biot-Savart operator, and use it to help us obtain formula (2) of Theorem 3 for the Green's operator. Note that this is the reverse of the order of derivation that we used in left-translation format.

It is interesting to contrast the two formats for thinking about and calculating with vector fields on $S^3$.

Left-translation format takes advantage of the group structure on $S^3$ and has the following very attractive features:

• If we start with a tangent vector $V$ at a single point $x \in S^3$ and left-translate it to each possible point $y \in S^3$, we get a vector field $W(y) = L_{yx^{-1}} V$ which is left-invariant, and tangent to the fibres of a Hopf fibration of $S^3$ by parallel great circles.

• This left-invariant vector field $W(y) = L_{yx^{-1}} V$ is divergence-free, and is a curl-eigenfield with eigenvalue $-2$.

We have already seen, in the study of the calculus of vector convolutions in $S^3$, how "computationally convenient" these two features are.

But left-translation format, however attractive and convenient, is specific to $S^3$ because of its group structure. Parallel transport, by contrast, makes sense on any Riemannian manifold, although it is limited by the non-uniqueness of geodesic connections between distant points.

Consider the following features of parallel transport on $S^3$ :

• If we start with a tangent vector $V$ at a single point $x \in S^3$ and then parallel transport it to each possible point $y \in S^3$, we get a vector field $W(y) = P_{yx} V$ which is defined for all $y \in S^3$ *except* for the point $y = -x$ antipodal to $x$.

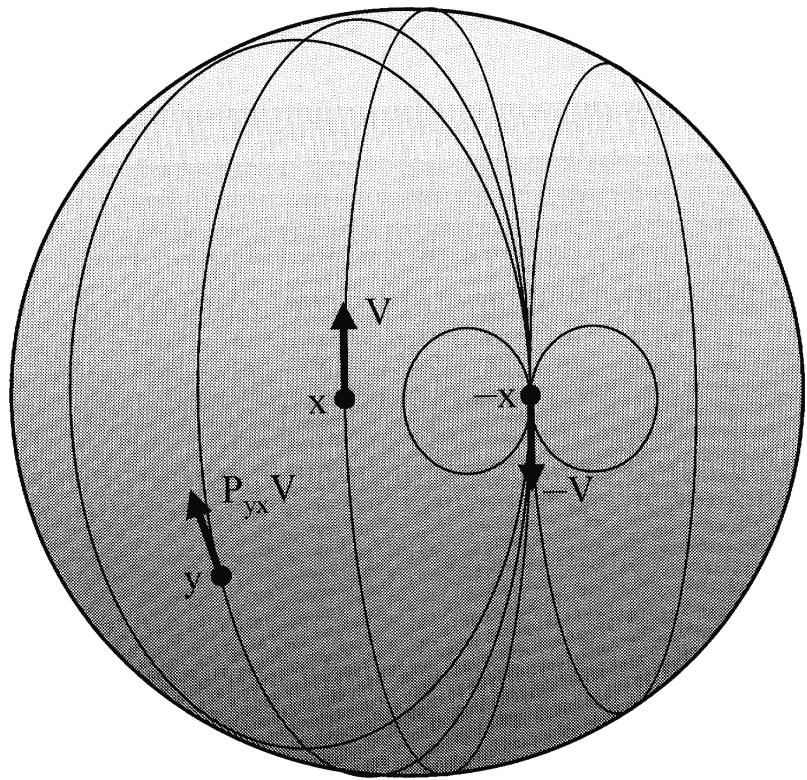

**The orbits of the vector field $P_{yx}V$ are all the oriented circles, great and small, which are tangent at the point $-x$ to the vector $-V$.**

• The vector field $P_{yx} V$ is neither divergence-free nor a curl eigenfield. We will obtain and use explicit formulas for its divergence and curl.

# 14. Geodesics in $S^3$.

If $x, y \in R^4$, we will write $<x, y>$ for the standard Euclidean inner product

$$<x, y> = x_0 y_0 + x_1 y_1 + x_2 y_2 + x_3 y_3$$

in $R^4$. Then

$$S^3 = \{x \in R^4 : <x, x> = 1\}.$$

We will freely identify points in $R^4$ with the vectors at the origin that point to them, and we will also identify the tangent spaces at various points in $R^4$ in the usual way, so that we can say $v \in T_x S^3$ if and only if $<x, v> = 0$. We will always write $<v, w>$ for the inner product of vectors when we're thinking of them as belonging to $R^4$, and reserve the notation $v \bullet w$ for the induced inner product on $T_x S^3$.

Suppose $x \in S^3$ so that $<x, x> = 1$, and suppose $v$ is a unit vector tangent to $S^3$ at $x$, so that $<x, v> = 0$ and $<v, v> = 1$. Then the great circle

$$G(t) = \cos t \, x + \sin t \, v$$

is a geodesic in $S^3$ parametrized by arc length.

Since $<x, G(t)> = \cos t$, we have that for any pair of points $x, y \in S^3$, their inner product $<x, y>$ is equal to the cosine of the distance (along $S^3$) from $x$ to $y$.

Now let x and y be any pair of distinct, non-antipodal points in $S^3$. Then the vector $v = y - <x, y> x$ is nonzero and orthogonal to x, and $|v|^2 = <v, v> = 1 - <x, y>^2$. Therefore the unique, length-minimizing geodesic from x to y is given by

$$G(t) = \cos t\ x + \sin t\ (v / |v|).$$

If $\alpha$ is the distance from x to y, so that $\cos \alpha = <x, y>$, then we can rewrite this as

$$G(t) = \cos t\ x + \sin t\ ((y - \cos \alpha\ x) / \sin \alpha).$$

When $t = \alpha$, we have $G(\alpha) = y$ and

$$G'(\alpha) = -\sin \alpha\ x + \cos \alpha\ ((y - \cos \alpha\ x) / \sin \alpha)$$

$$= (\cos \alpha\ y - x) / \sin \alpha.$$

Because G is parametrized by arc-length and is length-minimizing, and $\alpha$ is the distance in $S^3$ from x to y, we deduce that

$$\nabla_y \alpha(x, y) = G'(\alpha) = (\cos \alpha\ y - x) / \sin \alpha.$$

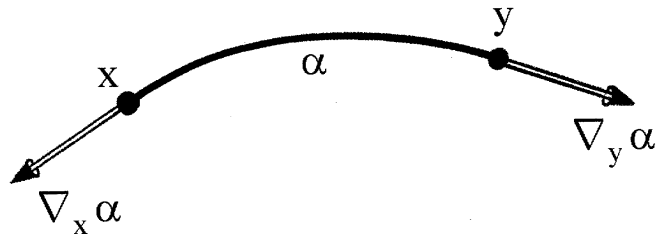

**The gradient vectors of $\alpha(x, y)$ at x and at y**

We note that

$$P_{xy} \nabla_y \alpha(x, y) = -\nabla_x \alpha(x, y).$$

# 15. The vector triple product in $R^4$.

Let $X_0 = (1, 0, 0, 0)$, $X_1 = (0, 1, 0, 0)$, $X_2 = (0, 0, 1, 0)$ and $X_3 = (0, 0, 0, 1)$ be the standard orthonormal basis for $R^4$. Then, for $x, v, w \in R^4$, we define the *vector triple product*

$$[x, v, w] = \det \begin{pmatrix} X_0 & X_1 & X_2 & X_3 \\ v_0 & v_1 & v_2 & v_3 \\ w_0 & w_1 & w_2 & w_3 \\ X_0 & X_1 & X_2 & X_3 \end{pmatrix} \in R^4.$$

Note that the first three rows of the above "determinant" consist of real numbers, while the last row consists of vectors in $R^4$. The value $[x, v, w]$ of this determinant is a vector, orthogonal to $x$, $v$ and $w$, whose length is equal to the volume of the 3-dimensional parallelepiped spanned by $x$, $v$ and $w$. Note that $[x, v, w]$ is an alternating, multilinear function of its three arguments.

For $y \in R^4$, we have

$$\langle [x, v, w], y \rangle = \det \begin{pmatrix} x_0 & x_1 & x_2 & x_3 \\ v_0 & v_1 & v_2 & v_3 \\ w_0 & w_1 & w_2 & w_3 \\ y_0 & y_1 & y_2 & y_3 \end{pmatrix}.$$

If $A \in SO(4)$, then

$$[Ax, Av, Aw] = A[x, v, w].$$

# 16. Explicit formula for parallel transport in $S^3$.

Let $x$ and $y$ be any pair of distinct, non-antipodal points in $S^3 \subset R^4$, equivalently, any pair of linearly independent unit vectors in $R^4$. Then there is a unique element $M \in SO(4)$ that maps $x$ to $y$ and leaves fixed the two-dimensional subspace of $R^4$ consisting of all vectors orthogonal to both $x$ and $y$. For $v \in R^4$, this mapping is given by

$$M(v) = v - \frac{<x+y, v>}{1 + <x, y>} x + \left\langle \frac{1 + 2<x, y>}{1 + <x, y>} x - \frac{1}{1 + <x, y>} y, v \right\rangle y.$$

Direct calculation shows that $M \in SO(4)$, and is therefore the desired isometry.

If $v \in T_x S^3$, then $<v, x> = 0$ and the formula for $M(v)$ simplifies to

$$M(v) = v - \frac{<y, v>}{1 + <x, y>} (x + y).$$

In this case, $M(v)$ is the result of parallel transport of $v$ from $x$ to $y$, so we will write

$$P_{yx} v = M(v) = v - \frac{<y, v>}{1 + <x, y>} (x + y).$$

# 17. The method of moving frames in $S^3$.

The computations which lead to the various integral formulas on $S^3$ in parallel transport format will be carried out using "moving frames" which are obtained via parallel transport from a single orthonormal frame, as follows.

Let $x = (1, 0, 0, 0)$, and let $X_1$, $X_2$, $X_3$ be the usual orthonormal basis for the tangent space $T_xS^3$ to $S^3$ at $x$.

For each point $y = (y_0, y_1, y_2, y_3) \neq -x$ in $S^3$, we parallel transport this orthonormal basis for $T_xS^3$ along the unique shortest geodesic (great circle arc) from $x$ to $y$ to obtain an orthonormal basis

$$E_1 = P_{yx}X_1, \quad E_2 = P_{yx}X_2, \quad E_3 = P_{yx}X_3$$

for $T_yS^3$.

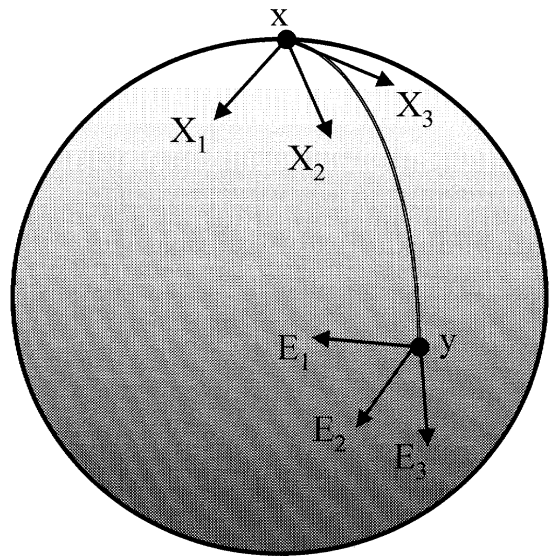

**Moving frames**

Using the formula for $P_{yx}V$ from the previous section, we calculate that

$$E_1 = (-y_1,\ 1 - y_1^2/(1+y_0),\ -y_1y_2/(1+y_0),\ -y_1y_3/(1+y_0))$$

$$E_2 = (-y_2,\ -y_1y_2/(1+y_0),\ 1 - y_2^2/(1+y_0),\ -y_2y_3/(1+y_0))$$

$$E_3 = (-y_3,\ -y_1y_3/(1+y_0),\ -y_2y_3/(1+y_0),\ 1 - y_3^2/(1+y_0)).$$

Next, let $\theta_1$, $\theta_2$, $\theta_3$ denote the 1-forms on $S^3 - \{-x\}$ which are dual to $E_1$, $E_2$, $E_3$, in the sense that $\theta_i(V) = \langle E_i, V \rangle$ for all vectors $V$ and for $i = 1, 2, 3$.

We can write

$$\theta_1 = -y_1\,dy_0 + (1 - y_1^2/(1+y_0))\,dy_1 - y_1y_2/(1+y_0)\,dy_2 - y_1y_3/(1+y_0)\,dy_3$$

$$\theta_2 = -y_2\,dy_0 - y_1y_2/(1+y_0)\,dy_1 + (1 - y_2^2/(1+y_0))\,dy_2 - y_2y_3/(1+y_0)\,dy_3$$

$$\theta_3 = -y_3\,dy_0 - y_1y_3/(1+y_0)\,dy_1 - y_2y_3/(1+y_0)\,dy_2 + (1 - y_3^2/(1+y_0))\,dy_3,$$

where we understand the right hand sides as restricted to $S^3$.

Since $y_0^2 + y_1^2 + y_2^2 + y_3^2 = 1$ on $S^3$, we can add a fourth equation to the three above,

$$0 = y_0\,dy_0 + y_1\,dy_1 + y_2\,dy_2 + y_3\,dy_3.$$

Using this fourth equation, we can simplify the first three to read

$$\theta_1 = -y_1/(1+y_0)\, dy_0 + dy_1$$

$$\theta_2 = -y_2/(1+y_0)\, dy_0 + dy_2$$

$$\theta_3 = -y_3/(1+y_0)\, dy_0 + dy_3\,.$$

We then invert this system to obtain

$$dy_0 = -y_1\,\theta_1 - y_2\,\theta_2 - y_3\,\theta_3$$

$$dy_1 = \theta_1 - (y_1/(1+y_0))\,(y_1\,\theta_1 + y_2\,\theta_2 + y_3\,\theta_3)$$

$$dy_2 = \theta_2 - (y_2/(1+y_0))\,(y_1\,\theta_1 + y_2\,\theta_2 + y_3\,\theta_3)$$

$$dy_3 = \theta_3 - (y_3/(1+y_0))\,(y_1\,\theta_1 + y_2\,\theta_2 + y_3\,\theta_3)\,.$$

With these equations in hand, we easily obtain

$$d\theta_1 = y_2/(1+y_0)\,\theta_1 \wedge \theta_2 - y_3/(1+y_0)\,\theta_3 \wedge \theta_1$$

$$d\theta_2 = y_3/(1+y_0)\,\theta_2 \wedge \theta_3 - y_1/(1+y_0)\,\theta_1 \wedge \theta_2$$

$$d\theta_3 = y_1/(1+y_0)\,\theta_3 \wedge \theta_1 - y_2/(1+y_0)\,\theta_2 \wedge \theta_3\,.$$

## 18. Curls and divergences.

If x is a fixed point on $S^3$ and V is a tangent vector to $S^3$ at x, then the parallel transport $P_{yx}V$ of V to the point y is defined for all $y \neq -x$, and may be viewed as a vector field on $S^3 - \{-x\}$, as illustrated in section 13.

We need to know the curl and divergence of this vector field.

**PROPOSITION 18.1.**
$$\nabla_y \times P_{yx}V = \frac{[y, x, V]}{1 + <x, y>} .$$

**PROPOSITION 18.2.**
$$\nabla_y \bullet P_{yx}V = \frac{-2 <y, V>}{1 + <x, y>} .$$

To take curls and divergences of vector fields on an oriented Riemannian 3-manifold, we make use of the duality between vector fields and forms as follows.

If U is a vector field, we can convert it to the dual 1-form $\downarrow U$ which is defined by $(\downarrow U)(W) = <U, W>$ for all vector fields W. If µ is a 1-form, we can convert it to the dual vector field $\uparrow \mu$ which is defined by $<\uparrow \mu, W> = \mu(W)$ for all vector fields W. The down and up arrows remind us of lowering and raising indices in tensorial notation.

If U is a vector field, we can convert it to the dual 2-form $\Downarrow U$ which is defined by $(\Downarrow U)(Y, Z) = \text{vol}(U, Y, Z)$ for all vector fields Y and Z. If µ is a 2-form, we can convert it to the dual vector field $\Uparrow \mu$ which is defined by $\text{vol}(\Uparrow \mu, Y, Z) = \mu(Y, Z)$ for all vector fields Y and Z.

For example, $\downarrow E_i = \theta_i$ and $\uparrow \theta_i = E_i$ for i = 1, 2, 3. Likewise $\Downarrow E_1 = \theta_2 \wedge \theta_3$, $\Downarrow E_2 = \theta_3 \wedge \theta_1$ and $\Downarrow E_3 = \theta_1 \wedge \theta_2$, and $\Uparrow$ is simply the inverse of this.

With these notations set, the curl of a vector field $W$ is defined, as usual, by

$$\nabla \times W = \Uparrow d \downarrow W,$$

that is, we convert $W$ to the 1-form $\downarrow W$, take its exterior derivative to get the 2-form $d \downarrow W$, and then convert this to the vector field $\Uparrow d \downarrow W$.

For example, suppose $W$ is the vector field $E_1 = P_{yx}X_1$. Then

$$\nabla_y \times E_1 = \Uparrow d \downarrow E_1 = \Uparrow d\, \theta_1$$

$$= \Uparrow (y_2/(1 + y_0)\, \theta_1 \wedge \theta_2 - y_3/(1 + y_0)\, \theta_3 \wedge \theta_1)$$

$$= y_2/(1 + y_0)\, E_3 - y_3/(1 + y_0)\, E_2.$$

And likewise

$$\nabla_y \times E_2 = y_3/(1 + y_0)\, E_1 - y_1/(1 + y_0)\, E_3$$

$$\nabla_y \times E_3 = y_1/(1 + y_0)\, E_2 - y_2/(1 + y_0)\, E_1.$$

These formulas can be viewed as special instances of the formula given above in Proposition 18.1. For example, if $y = (y_0, y_1, y_2, y_3)$, $x = (1, 0, 0, 0)$ and $V = X_1$, then one easily calculates that

$$[y, x, V] = y_2\, E_3 - y_3\, E_2.$$

Dividing through by $1 + \langle x, y \rangle = 1 + y_0$ gives the formula for $\nabla_y \times E_1$ appearing above, and so verifies Proposition 18.1 in this instance.

Both sides of the formula in Proposition 18.1 are linear in $V$, so we can certainly rescale $V$ to be a unit vector. Also, both sides of this formula are invariant under an orientation preserving isometry of $S^3$, so we can move $x$ to the point $(1, 0, 0, 0)$ and $V$ to the tangent vector $E_1$.

Thus the above verification of Proposition 18.1 in this single instance also serves as its proof.

The divergence of a vector field $W$ is defined by

$$\nabla \bullet W = * d \Downarrow W,$$

where $*$ is the Hodge star operator which converts the 3-form $f\, d(vol)$ to the function $f$. Thus the divergence of the vector field $W$ is calculated by converting it to the 2-form $\Downarrow W$, taking its exterior derivative to get the 3-form $d \Downarrow W$, and then converting this to the function $* d \Downarrow W$.

For example, suppose $W$ is the vector field $E_1 = P_{yx} X_1$. Then

$$\nabla \bullet E_1 = * d \Downarrow E_1 = * d\, (\theta_2 \wedge \theta_3)$$

$$= * (d\theta_2 \wedge \theta_3 - \theta_2 \wedge d\theta_3)$$

$$= * (-2\, y_1/(1 + y_0))\, \theta_1 \wedge \theta_2 \wedge \theta_3$$

$$= -2\, y_1/(1 + y_0),$$

using the formulas for $d\theta_2$ and $d\theta_3$ appearing in the previous section.

The above example can be viewed as a special instance of the formula given in Proposition 18.2, since

$$<y, V> \; = \; <y, X_1> \; = \; y_1 \quad \text{and} \quad 1 + <x, y> \; = \; 1 + y_0.$$

And just as in the case of the curl formula of Proposition 18.1, this single verification of Proposition 18.2 also serves as its proof.

## 19. Statement of the Key Lemma.

The following result plays a lead role in obtaining the formulas for the Biot-Savart and Green's operators on $S^3$ in parallel transport format.

**KEY LEMMA, spherical version.**

$$\nabla_y \times \{P_{yx} V(x) \times \nabla_y \varphi\} - \nabla_y \{V(x) \bullet \nabla_x (\cos \alpha \ \varphi)\}$$

$$= (\Delta \varphi - \varphi) (V(x) - <V(x), y> y) .$$

This formula is quite flexible: it is useful whether $V$ is a smooth vector field defined on all of $S^3$, or just on a bounded subdomain $\Omega$ of $S^3$, or just along a smooth closed curve $K$ in $S^3$, or even just at the single point $x$.

As usual, $\alpha$ is the distance on $S^3$ between $x$ and $y$. The kernel function $\varphi(\alpha)$ is any smooth function, which typically blows up at $\alpha = 0$ and is asymptotic there to $-1/(4\pi\alpha)$.

The Key Lemma can be understood intuitively in the special case that $\varphi(\alpha)$ is the fundamental solution of the operator $\varphi \to \Delta\varphi - \varphi$ on $S^3$, in which case we will see that $\cos \alpha \ \varphi(\alpha)$ is the fundamental solution of the scalar Laplacian there.

With this choice of kernel function, the quantity $P_{yx} V(x) \times \nabla_y \varphi$ can be viewed as the contribution made by a little bit of current $V(x)$ at the location $x$ to the magnetic field at the location $y$. When integrated over $S^3$ with respect to $x$, the result will be the magnetic field $B$ at $y$. Thus the first term on the left hand side of the Key Lemma integrates to $\nabla \times B$.

In similar fashion, the term $\nabla_y \{V(x) \bullet \nabla_x (\cos \alpha \ \varphi)\}$ integrates to the time derivative $\partial E/\partial t$ of the electric field at $y$ due to the accumulation or dispersion of charge by the current flow $V$, while the right hand side of the Key Lemma integrates to the current flow $V$ itself.

Thus the Key Lemma, for this special choice of kernel function $\varphi$, can be viewed as a pointwise or infinitesimal version of Maxwell's equation

$$\nabla \times B - \partial E/\partial t \ = V .$$

## 20. The first term on the LHS of the Key Lemma.

The plan for the proof of the Key Lemma is simple. Applying an appropriate rotation of $S^3$, we can assume that $x$ is the point $(1, 0, 0, 0)$ and that, after scaling, $V(x) = X_1 = (0, 1, 0, 0)$. At each point $y = (y_0, y_1, y_2, y_3)$ of $S^3 - \{-x\}$, we will expand both sides of the Key Lemma in terms of the orthonormal basis $E_1$, $E_2$ and $E_3$ for $T_y S^3$, and prove equality component-wise.

We begin with the first term on the LHS of the equation, and note that

$$P_{yx} V(x) \times \nabla_y \varphi = P_{yx} V(x) \times (-P_{yx} \nabla_x \varphi) = P_{yx} (-V(x) \times \nabla_x \varphi).$$

From section 14, we have that

$$\nabla_x \alpha = (\cos \alpha \ x - y) / \sin \alpha$$

$$= ((y_0, 0, 0, 0) - (y_0, y_1, y_2, y_3)) / \sin \alpha$$

$$= (-1 / \sin \alpha)(y_1 X_1 + y_2 X_2 + y_3 X_3).$$

Hence

$$\nabla_x \varphi = \varphi'(\alpha) \nabla_x \alpha = (-\varphi'(\alpha) / \sin \alpha)(y_1 X_1 + y_2 X_2 + y_3 X_3).$$

Thus

$$-V(x) \times \nabla_x \varphi = (-X_1) \times (-\varphi'(\alpha) / \sin \alpha)(y_1 X_1 + y_2 X_2 + y_3 X_3)$$

$$= (\varphi'(\alpha) / \sin \alpha)(-y_3 X_2 + y_2 X_3).$$

Hence

$$P_{yx} V(x) \times \nabla_y \varphi = P_{yx} (-V(x) \times \nabla_x \varphi)$$

$$= (\varphi'(\alpha) / \sin \alpha)(-y_3 E_2 + y_2 E_3).$$

Now we take curls and evaluate the first term on the LHS of the Key Lemma,

$$\nabla_y \times \{P_{yx} V(x) \times \nabla_y \varphi\} = \nabla_y \times \{(\varphi'(\alpha) / \sin \alpha)(-y_3 E_2 + y_2 E_3)\}.$$

We begin by computing the curl of the vector field

$$W = W(y) = -y_3 E_2 + y_2 E_3.$$

Following the recipe given in section 18, we have

$$\nabla_y \times W = \Uparrow d \Downarrow W.$$

Going part way, we have

$$d \Downarrow W = d(-y_3 \theta_2 + y_2 \theta_3)$$

$$= -dy_3 \wedge \theta_2 - y_3 d\theta_2 + dy_2 \wedge \theta_3 + y_2 d\theta_3.$$

We use the formulas from section 17 to express $dy_3$ and $dy_2$ in terms of $\theta_1$, $\theta_2$ and $\theta_3$, and compute that

$$d \Downarrow W = 2(y_0 + y_1^2/(1+y_0))\theta_2 \wedge \theta_3 + 2(y_1 y_2 / (1+y_0))\theta_3 \wedge \theta_1$$

$$+ 2(y_1 y_3 / (1+y_0))\theta_1 \wedge \theta_2.$$

Hence

$$\nabla_y \times W = \Uparrow d \Downarrow W = 2(y_0 + y_1^2/(1+y_0)) E_1 + 2(y_1 y_2 / (1+y_0)) E_2$$

$$+ 2(y_1 y_3 / (1+y_0)) E_3.$$

Now to compute $\nabla_y \times \{(\varphi'(\alpha) / \sin \alpha)(-y_3 E_2 + y_2 E_3)\}$, we use the formula

$$\nabla \times (fW) = f(\nabla \times W) + \nabla f \times W,$$

with $f(\alpha) = \varphi'(\alpha) / \sin \alpha$.

We find that

$$\nabla_y \times \{P_{yx}V(x) \times \nabla_y \varphi\} = \nabla_y \times \{(\varphi'(\alpha)/\sin\alpha)(-y_3 E_2 + y_2 E_3)\}$$

$$= \{(\varphi'/\sin\alpha)\, 2\,(y_0 + y_1^2/(1+y_0)) + (\varphi''/\sin^2\alpha - \varphi'\cos\alpha/\sin^3\alpha)(y_2^2 + y_3^2)\}\, E_1$$

$$+ \{(\varphi'/\sin\alpha)\, 2\, y_1 y_2/(1+y_0) + (\varphi''/\sin^2\alpha - \varphi'\cos\alpha/\sin^3\alpha)(-y_1 y_2)\}\, E_2$$

$$+ \{(\varphi'/\sin\alpha)\, 2\, y_1 y_3/(1+y_0) + (\varphi''/\sin^2\alpha - \varphi'\cos\alpha/\sin^3\alpha)(-y_1 y_3)\}\, E_3.$$

This puts the first term on the LHS of the Key Lemma in the form we want it.

## 21. The second term on the LHS of the Key Lemma.

Now, just as we did for the first term, we want to express the second term

$$-\nabla_y \{V(x) \bullet \nabla_x (\cos \alpha \ \varphi)\}$$

as a linear combination of the orthonormal basis vectors $E_1$, $E_2$ and $E_3$ for $T_y S^3$. We keep in mind that $x = (1, 0, 0, 0)$ and that $V(x) = X_1 = (0, 1, 0, 0)$.

To begin,

$$\nabla_x (\cos \alpha \ \varphi) = (\cos \alpha \ \varphi)' \nabla_x \alpha$$

$$= (\cos \alpha \ \varphi)' (-1/\sin \alpha) (y_1 X_1 + y_2 X_2 + y_3 X_3).$$

Then taking the dot product of this with $V(x) = X_1$ yields

$$V(x) \bullet \nabla_x (\cos \alpha \ \varphi) = (\cos \alpha \ \varphi)' (-1/\sin \alpha) y_1$$

$$= (-\varphi' \cos \alpha / \sin \alpha + \varphi) y_1.$$

Thus

$$-\nabla_y \{V(x) \bullet \nabla_x (\cos \alpha \ \varphi)\} = \nabla_y \{(\varphi' \cos \alpha / \sin \alpha - \varphi)\} y_1.$$

To calculate this, recall that

$$\nabla_y \alpha = (1/\sin \alpha) (y_1 E_1 + y_2 E_2 + y_3 E_3),$$

and note that

$$\nabla_y y_1 = \uparrow dy_1 = \uparrow \{\theta_1 - (y_1/(1 + y_0)) (y_1 \theta_1 + y_2 \theta_2 + y_3 \theta_3)\}$$

$$= E_1 - (y_1/(1 + y_0)) (y_1 E_1 + y_2 E_2 + y_3 E_3).$$

We then calculate the second term on the LHS of the Key Lemma:

$-\nabla_y \{V(x) \bullet \nabla_x (\cos \alpha \; \varphi)\}$

$$= \{ (\varphi'' \cos \alpha / \sin^2 \alpha - \varphi'(1/\sin \alpha + 1/\sin^3 \alpha)) y_1^2$$
$$+ (\varphi' \cos \alpha / \sin \alpha - \varphi)(1 - y_1^2/(1+y_0)) \} E_1$$
$$+ \{ (\varphi'' \cos \alpha / \sin^2 \alpha - \varphi'(1/\sin \alpha + 1/\sin^3 \alpha)) y_1 y_2$$
$$- (\varphi' \cos \alpha / \sin \alpha - \varphi) y_1 y_2/(1+y_0) \} E_2$$
$$+ \{ (\varphi'' \cos \alpha / \sin^2 \alpha - \varphi'(1/\sin \alpha + 1/\sin^3 \alpha)) y_1 y_3$$
$$- (\varphi' \cos \alpha / \sin \alpha - \varphi) y_1 y_3/(1+y_0) \} E_3 .$$

If we combine the two terms on the LHS of the Key Lemma and simplify, we get

$\nabla_y \times \{P_{yx} V(x) \times \nabla_y \varphi\} - \nabla_y \{V(x) \bullet \nabla_x (\cos \alpha \; \varphi)\}$

$$= \{(\varphi'' + 2\varphi' \cos \alpha / \sin \alpha - \varphi)(1 - y_1^2/(1+y_0)) \} E_1$$
$$+ \{(\varphi'' + 2\varphi' \cos \alpha / \sin \alpha - \varphi)(-y_1 y_2/(1+y_0)) \} E_2$$
$$+ \{(\varphi'' + 2\varphi' \cos \alpha / \sin \alpha - \varphi)(-y_1 y_3/(1+y_0)) \} E_3 .$$

## 22. The RHS of the Key Lemma.

Now we express the RHS

$$(\Delta\varphi - \varphi)(V(x) - <V(x), y> y)$$

as a linear combination of $E_1$, $E_2$ and $E_3$.

Recall that

$$\Delta\varphi = \varphi'' + 2\varphi' \cos\alpha/\sin\alpha,$$

and hence

$$\Delta\varphi - \varphi = \varphi'' + 2\varphi' \cos\alpha/\sin\alpha - \varphi.$$

Next, note that $V(x) = (0, 1, 0, 0)$ can be viewed as the gradient in $R^4$ of the coordinate function $y_1$, and that $V(x) - <V(x), y> y$ is the orthogonal projection of this onto the tangent plane $T_y S^3$ to $S^3$ at $y$, and is therefore the gradient on $S^3$ of $y_1$. Thus

$$V(x) - <V(x), y> y = \nabla_y y_1$$

$$= E_1 - (y_1/(1 + y_0))(y_1 E_1 + y_2 E_2 + y_3 E_3),$$

according to our calculation in the preceding section.

Hence the RHS of the Key Lemma is given by

$$(\Delta\varphi - \varphi)(V(x) - <V(x), y> y)$$

$$= \{(\varphi'' + 2\varphi' \cos\alpha/\sin\alpha - \varphi)(1 - y_1^2/(1 + y_0))\} E_1$$

$$+ \{(\varphi'' + 2\varphi' \cos\alpha/\sin\alpha - \varphi)(-y_1 y_2/(1 + y_0))\} E_2$$

$$+ \{(\varphi'' + 2\varphi' \cos\alpha/\sin\alpha - \varphi)(-y_1 y_3/(1 + y_0))\} E_3.$$

This matches the LHS, as given in the previous section, and so completes the proof of the Key Lemma.

# 23. Maxwell's equation on $S^3$.

We now apply the Key Lemma to the proof of Maxwell's equation

$$\nabla \times B = V + \partial E/\partial t$$

on $S^3$. In the Key Lemma, we use the kernel function

$$\varphi(\alpha) = (-1/4\pi^2)(\pi - \alpha) \csc \alpha ,$$

which, as we noted earlier, is a fundamental solution of the shifted Laplacian on $S^3$,

$$\Delta\varphi - \varphi = \delta .$$

With this choice of kernel function, we have that

$$(\cos \alpha)\, \varphi(\alpha) = \varphi_0(\alpha) = (-1/4\pi^2)(\pi - \alpha) \cot \alpha ,$$

which is the fundamental solution of the Laplacian on $S^3$,

$$\Delta\varphi_0 = \delta - 1/(2\pi^2) .$$

With these choices, the Key Lemma reads

$$\nabla_y \times \{P_{yx}V(x) \times \nabla_y\varphi\} - \nabla_y \{V(x) \bullet \nabla_x \varphi_0\}$$

$$= \delta(x, y)\,(V(x) - <V(x), y> y) .$$

Keeping these choices of kernel functions, we now define the Biot-Savart operator by formula (2) of Theorem 2,

$$BS(V)(y) = \int_{S^3} P_{yx}V(x) \times \nabla_y \varphi(x, y)\, dx .$$

Then we integrate the Key Lemma over $S^3$ and get

$$\int_{S^3} \nabla_y \times \{P_{yx}V(x) \times \nabla_y\varphi\}\, dx - \int_{S^3} \nabla_y \{V(x) \bullet \nabla_x \varphi_0\}\, dx$$

$$= \int_{S^3} \delta(x, y)\,(V(x) - <V(x), y> y)\, dx .$$

The integral on the right hand side above equals

$$V(y) - \langle V(y), y \rangle y = V(y),$$

since $V(y)$ is a tangent vector to $S^3$ at $y$, and hence orthogonal to $y$ in $R^4$.

Taking the curl and the gradient outside the two integrals on the left-hand side, we then get

$$\nabla_y \times \int_{S^3} P_{yx} V(x) \times \nabla_y \varphi(x, y) \, dx - \nabla_y \int_{S^3} V(x) \bullet \nabla_x \varphi_0(x, y) \, dx = V(y),$$

or

$$\nabla_y \times BS(V)(y) = V(y) + \nabla_y \int_{S^3} V(x) \bullet \nabla_x \varphi_0(x, y) \, dx.$$

As indicated in section 3, this is Maxwell's equation

$$\nabla \times B = V + \partial E/\partial t.$$

This explains our remark in section 19 that, for this special choice of kernel function $\varphi$, the Key Lemma can be viewed as a pointwise or infinitesimal version of Maxwell's equation.

## 24. Proof of Theorem 2, formula (2).

Let $V$ be a smooth vector field on $S^3$. Thinking of $V$ as a steady current flow, we defined the corresponding magnetic field $BS(V)$ by

$$BS(V)(y) = \int_{S^3} P_{yx} V(x) \times \nabla_y \varphi(x, y) \, dx,$$

where $\varphi(\alpha) = (-1/4\pi^2)(\pi - \alpha) \csc \alpha$.

We saw in the preceding section that, as a consequence of the Key Lemma, the field $BS(V)$ satisfies Maxwell's equation

$$\nabla_y \times BS(V)(y) = V(y) + \nabla_y \int_{S^3} V(x) \bullet \nabla_x \varphi_0(x, y) \, dx.$$

To complete the proof of Theorem 2, formula (2), we need only show that $BS(V)$ is divergence-free, and we do this as follows.

$$\nabla_y \bullet BS(V)(y) = \nabla_y \bullet \int_{S^3} P_{yx} V(x) \times \nabla_y \varphi(x, y) \, dx$$

$$= \int_{S^3} \nabla_y \bullet \{P_{yx} V(x) \times \nabla_y \varphi(x, y)\} \, dx,$$

and we will show that the integrand is identically zero. We use the standard formula

$$\nabla \bullet (A \times B) = (\nabla \times A) \bullet B - A \bullet (\nabla \times B),$$

and also the formula from Proposition 18.1,

$$\nabla_y \times P_{yx} V(x) = [y, x, V] / (1 + \langle x, y \rangle).$$

Then

$$\nabla_y \bullet \{P_{yx}V(x) \times \nabla_y\varphi(x, y)\}$$

$$= \{\nabla_y \times P_{yx}V(x)\} \bullet \nabla_y\varphi(x, y) - P_{yx}V(x) \bullet \nabla_y \times \nabla_y\varphi(x, y)$$

$$= \{\nabla_y \times P_{yx}V(x)\} \bullet \nabla_y\varphi(x, y)$$

$$= \{[y, x, V(x)] / (1 + <x, y>)\} \bullet \nabla_y\varphi(x, y)$$

$$= 0,$$

since $[y, x, V(x)]$ is orthogonal to the plane spanned by $x$ and $y$, while $\nabla_y\varphi(x, y)$ lies in that plane.

Thus the proposed magnetic field $BS(V)$ is indeed divergence-free.

This completes the proof of Theorem 2, formula (2).

We note that the above argument for the divergence-free character of $BS(V)$ does not depend on the particular choice of kernel function $\varphi$.

# 25. The Green's operator on $S^3$ in parallel transport format.

We turn now to the derivation of formula (2) of Theorem 3 for the Green's operator on $S^3$ in parallel transport format,

$$Gr(V)(y) = \int_{S^3} P_{yx}V(x)\, \varphi_2(x, y)\, dx + \nabla_y \int_{S^3} P_{yx}V(x) \bullet \nabla_y \varphi_3(x, y)\, dx ,$$

where

$$\varphi_2(\alpha) = (-1/4\pi^2)(\pi - \alpha) \csc \alpha + (1/8\pi^2)(\pi - \alpha)^2 / (1 + \cos \alpha)$$

and

$$\varphi_3(\alpha) = (-1/24)(\pi - \alpha) \cot \alpha - (1/16\pi^2)\, \alpha\, (2\pi - \alpha)$$

$$+ (1/8\pi^2) \int_\alpha^\pi \left((\pi - \alpha)^3 / (3 \sin^2\alpha)\right) + \left((\pi - \alpha)^2 / \sin \alpha\right)\, d\alpha .$$

The naming of the kernel functions $\varphi_2$ and $\varphi_3$ above leaves room in the $\varphi$-family for the fundamental solution of the Laplacian on $S^3$,

$$\varphi_0(\alpha) = (-1/4\pi^2)(\pi - \alpha) \cot \alpha \qquad \Delta\varphi_0 = \delta - 1/2\pi^2$$

and for the kernel function of the Biot-Savart operator in parallel transport format,

$$\varphi_1(\alpha) = (-1/4\pi^2)(\pi - \alpha) \csc \alpha \qquad \Delta\varphi_1 - \varphi_1 = \delta .$$

## 26. Plan of the proof.

Consider the usual vector convolutions, this time in parallel transport format:

$$A(V, \varphi)(y) = \int_{S^3} P_{yx} V(x) \, \varphi(x, y) \, dx$$

$$B(V, \varphi)(y) = \int_{S^3} P_{yx} V(x) \times \nabla_y \varphi(x, y) \, dx$$

$$g(V, \varphi)(y) = \int_{S^3} P_{yx} V(x) \bullet \nabla_y \varphi(x, y) \, dx$$

$$G(V, \varphi)(y) = \nabla_y \, g(V, \varphi)(y) = \nabla_y \int_{S^3} P_{yx} V(x) \bullet \nabla_y \varphi(x, y) \, dx \, .$$

*Step 1.* First we will determine the unknown kernel function $\varphi$ so as to make

$$\nabla_y \times A(V, \varphi)(y) = -B(V, \varphi_1)(y) = -BS(V)(y) \, .$$

Insisting that $\varphi$ have a singularity only at $\alpha = 0$ and *not* at $\alpha = \pi$ leads to the unique solution $\varphi = \varphi_2$ given in the statement of the theorem above.

*Step 2.* Then, as a consequence of Maxwell's equation, we will have

$$-\nabla_y \times \nabla_y \times A(V, \varphi_2)(y) = \nabla_y \times BS(V)(y)$$

$$= V(y) - \nabla_y \int_{S^3} \nabla_x \bullet V(x) \, \varphi_0(x, y) \, dx \, .$$

Hence

$$\Delta_y \, A(V, \varphi_2)(y) = -\nabla_y \times \nabla_y \times A(V, \varphi_2)(y) + \nabla_y(\nabla_y \bullet A(V, \varphi_2)(y))$$

$$= V(y) - \nabla_y \int_{S^3} \nabla_x \bullet V(x) \, \varphi_0(x, y) \, dx + \nabla_y(\nabla_y \bullet A(V, \varphi_2)(y))$$

$$= V(y) + \text{some gradient field} \, .$$

We then find the kernel function $\varphi_3$ so that $\Delta_y \, G(V, \varphi_3)(y) = G(V, \Delta \varphi_3)(y)$ is the negative of the above gradient field, in which case defining

$$Gr(V)(y) = A(V, \varphi_2)(y) + G(V, \varphi_3)(y)$$

yields the desired equation $\Delta_y \, Gr(V)(y) = V(y) \, .$

## 27. Some formulas that we will use.

In the formulas below, $V$ is a single tangent vector at the single point $x \in S^3$.

(1) $$P_{yx}V = V - \frac{\langle y, V \rangle}{1 + \langle x, y \rangle}(x + y)$$

(2) $$\nabla_y \times P_{yx}V = \frac{[y, x, V]}{1 + \langle x, y \rangle}$$

(3) $$\nabla_y \bullet P_{yx}V = \frac{-2 \langle y, V \rangle}{1 + \langle x, y \rangle}$$

(4) $$P_{yx}V \times \nabla_y \varphi = [y, x, V]\, \varphi'(\alpha) / \sin \alpha .$$

Formula (1) appears in section 16, while formulas (2) and (3) appear in section 18.

Formula (4) in the special case that $x = (1, 0, 0, 0)$ and $V = X_1 = (0, 1, 0, 0)$ is obtained from the formula

$$P_{yx}V \times \nabla_y \varphi = (\varphi'(\alpha) / \sin \alpha)(-y_3 E_2 + y_2 E_3)$$

of section 20, together with the remark of section 18 that in this case we have

$$[y, x, V] = -y_3 E_2 + y_2 E_3 .$$

We leave it to the reader to confirm, as we did earlier in similar circumstances, that this single case of formula (4) can serve as its proof.

We see from formula (1) that $P_{yx}V$ lies in the 3-plane in $R^4$ spanned by $V$, $x$, $y$, and from formula (2) that its curl $\nabla_y \times P_{yx}V$ is orthogonal to this 3-plane, and hence orthogonal to $P_{yx}V$.

This is in marked contrast to the situation in left-translation format, where $\nabla_y \times L_{yx^{-1}}V = -2\, L_{yx^{-1}}V$ is parallel to $L_{yx^{-1}}V$.

We see from formulas (2) and (4) that $\nabla_y \times P_{yx}V$ and $P_{yx}V \times \nabla_y \varphi$ are parallel to one another. It is this fact which makes possible the solution of the equation

$$\nabla_y \times A(V, \varphi)(y) = -B(V, \varphi_1)(y) = -BS(V)(y)$$

for the unknown kernel function $\varphi$.

## 28. The curl of $A(V, \varphi)(y)$.

$$\nabla_y \times A(V, \varphi)(y) = \int_{S^3} \nabla_y \times \{P_{yx} V(x) \, \varphi(x, y)\} \, dx$$

$$= \int_{S^3} \{\nabla_y \times P_{yx} V(x)\} \, \varphi(x, y) \, - \, P_{yx} V(x) \times \nabla_y \varphi(x, y) \, dx$$

$$= \int_{S^3} \frac{[y, x, V(x)]}{1 + \langle x, y \rangle} \varphi(x, y) \, - \, [y, x, V(x)] \frac{\varphi'(\alpha)}{\sin \alpha} \, dx$$

$$= \int_{S^3} [y, x, V(x)] \left\{ \frac{\varphi}{1 + \cos \alpha} - \frac{\varphi'}{\sin \alpha} \right\} dx \, ,$$

thanks to formulas (2) and (4) above.

On the other hand,

$$BS(V)(y) = B(V, \varphi_1)(y) = \int_{S^3} P_{yx} V(x) \times \nabla_y \varphi_1(x, y) \, dx$$

$$= \int_{S^3} [y, x, V(x)] \frac{\varphi_1'(\alpha)}{\sin \alpha} \, dx \, ,$$

again by formula (4).

Thus, in order that

$$\nabla_y \times A(V, \varphi)(y) = -BS(V)(y) = -B(V, \varphi_1)(y) \, ,$$

we must choose $\varphi$ to satisfy the ODE

$$\frac{\varphi}{1 + \cos \alpha} - \frac{\varphi'}{\sin \alpha} = -\frac{\varphi_1'}{\sin \alpha} \, .$$

We rewrite this as

$$(1 + \cos \alpha) \varphi' - (\sin \alpha) \varphi = (1 + \cos \alpha) \varphi_1',$$

equivalently,

$$[(1 + \cos \alpha) \varphi]' = [(1 + \cos \alpha) \varphi_1]' + (\sin \alpha) \varphi_1$$

$$= [(1 + \cos \alpha) \varphi_1]' - (1/4\pi^2) (\pi - \alpha),$$

since $\varphi_1(\alpha) = (-1/4\pi^2) (\pi - \alpha) \csc \alpha$. Now we integrate both sides to get

$$(1 + \cos \alpha) \varphi = (1 + \cos \alpha) \varphi_1 + (1/8\pi^2) (\pi - \alpha)^2.$$

Dividing through by $(1 + \cos \alpha)$, we get

$$\varphi(\alpha) = \varphi_1(\alpha) + (1/8\pi^2) (\pi - \alpha)^2 / (1 + \cos \alpha)$$

$$= (-1/4\pi^2) (\pi - \alpha) \csc \alpha + (1/8\pi^2) (\pi - \alpha)^2 / (1 + \cos \alpha).$$

Note that $\varphi(\alpha)$ has a singularity only at $\alpha = 0$ and not at $\alpha = \pi$.

We call this solution $\varphi_2(\alpha)$, and with it we have completed Step 1:

$$\nabla_y \times A(V, \varphi_2)(y) = -B(V, \varphi_1)(y) = -BS(V)(y).$$

We take the negative curl of both sides to get

$$-\nabla_y \times \nabla_y \times A(V, \varphi_2)(y) = \nabla_y \times BS(V)(y)$$

$$= V(y) - \nabla_y \int_{S^3} (\nabla_x \bullet V(x)) \varphi_0(x, y) \, dx,$$

thanks to Maxwell's equation. Integrating by parts, we have

$$\int_{S^3} (\nabla_x \bullet V(x)) \varphi_0(x, y) \, dx = -\int_{S^3} V(x) \bullet \nabla_x \varphi_0(x, y) \, dx$$

$$= \int_{S^3} P_{yx} V(x) \bullet \nabla_y \varphi_0(x, y) \, dx$$

$$= g(V, \varphi_0)(y).$$

Hence

$$\Delta_y A(V, \varphi_2)(y) = -\nabla_y \times \nabla_y \times A(V, \varphi_2)(y) + \nabla_y(\nabla_y \bullet A(V, \varphi_2)(y))$$

$$= V(y) - \nabla_y \int_{S^3} \nabla_x \bullet V(x) \varphi_0(x, y) \, dx + \nabla_y(\nabla_y \bullet A(V, \varphi_2)(y))$$

$$= V(y) - G(V, \varphi_0)(y) + \nabla_y(\nabla_y \bullet A(V, \varphi_2)(y)).$$

Clearly the next step is to compute $\nabla_y \bullet A(V, \varphi_2)(y)$.

## 29. The divergence of $A(V, \varphi_2)$.

We begin as follows.

$$\nabla_y \bullet A(V, \varphi_2)(y) = \int_{S^3} \nabla_y \bullet \{P_{yx}V(x)\, \varphi_2(x, y)\} \, dx$$

$$= \int_{S^3} \{\nabla_y \bullet P_{yx}V(x)\}\, \varphi_2(x, y)\, dx + \int_{S^3} P_{yx}V(x) \bullet \nabla_y\varphi_2(x, y)\, dx.$$

The second integral on the RHS above is simply $g(V, \varphi_2)(y)$, so we put that aside.

Turning to the first integral and using formula (3) and the fact that $<x, V(x)> = 0$, we get

$$\int_{S^3} \{\nabla_y \bullet P_{yx}V(x)\}\, \varphi_2(x, y)\, dx = \int_{S^3} \frac{-2 <y, V(x)>}{1 + <x, y>} \varphi_2(x, y)\, dx$$

$$= \int_{S^3} <y, V(x)> \frac{-2\, \varphi_2(\alpha)}{1 + \cos\alpha}\, dx$$

$$= \int_{S^3} <\frac{\cos\alpha\, x - y}{\sin\alpha}, V(x)> \frac{2 \sin\alpha\, \varphi_2(\alpha)}{1 + \cos\alpha}\, dx$$

$$= \int_{S^3} <\nabla_x\alpha, V(x)> \frac{2 \sin\alpha\, \varphi_2(\alpha)}{1 + \cos\alpha}\, dx$$

$$= \int_{S^3} V(x) \bullet \frac{2 \sin\alpha\, \varphi_2(\alpha)}{1 + \cos\alpha}\, \nabla_x\alpha\, dx$$

$$= \int_{S^3} P_{yx}V(x) \bullet \frac{-2 \sin\alpha\, \varphi_2(\alpha)}{1 + \cos\alpha}\, \nabla_y\alpha\, dx$$

$$= \int_{S^3} P_{yx}V(x) \bullet \nabla_y\psi(\alpha)\, dx$$

$$= g(V, \psi)(y).$$

The new kernel function $\psi$ satisfies

$$\psi'(\alpha) = [-2 \sin \alpha / (1 + \cos \alpha)] \varphi_2(\alpha)$$

$$= [-2 \sin \alpha / (1 + \cos \alpha)] \left[ (-1/4\pi^2)(\pi - \alpha) \csc \alpha + (1/8\pi^2)(\pi - \alpha)^2 / (1 + \cos \alpha) \right]$$

$$= (1/2\pi^2)(\pi - \alpha)/(1 + \cos \alpha) - (1/4\pi^2)(\pi - \alpha)^2 \sin \alpha / (1 + \cos \alpha)^2,$$

and hence

$$\psi(\alpha) = (-1/4\pi^2)(\pi - \alpha)^2 / (1 + \cos \alpha).$$

Collecting information, we have

$$\nabla_y \bullet A(V, \varphi_2)(y)$$

$$= \int_{S^3} \{\nabla_y \bullet P_{yx} V(x)\} \varphi_2(x, y) \, dx + \int_{S^3} P_{yx} V(x) \bullet \nabla_y \varphi_2(x, y) \, dx$$

$$= g(V, \psi)(y) + g(V, \varphi_2)(y),$$

and hence

$$\nabla_y (\nabla_y \bullet A(V, \varphi_2)(y)) = G(V, \psi)(y) + G(V, \varphi_2)(y).$$

## 30. Step 2 - Finding and using the kernel function $\varphi_3$.

Continuing our calculation of the Laplacian of $A(V, \varphi_2)(y)$, we have

$$\Delta_y A(V, \varphi_2)(y) = V(y) - G(V, \varphi_0)(y) + \nabla_y(\nabla_y \bullet A(V, \varphi_2)(y))$$

$$= V(y) - G(V, \varphi_0)(y) + G(V, \psi)(y) + G(V, \varphi_2)(y)$$

$$= V(y) - G(V, \varphi_0 - \varphi_2 - \psi)(y).$$

We note that

$$\varphi_0 - \varphi_2 - \psi = (-1/4\pi^2)(\pi - \alpha)\cot\alpha$$

$$+ (1/4\pi^2)(\pi - \alpha)\csc\alpha - (1/8\pi^2)(\pi - \alpha)^2/(1 + \cos\alpha)$$

$$+ (1/4\pi^2)(\pi - \alpha)^2/(1 + \cos\alpha)$$

$$= (1/4\pi^2)(\pi - \alpha)(1 - \cos\alpha)/\sin\alpha + (1/8\pi^2)(\pi - \alpha)^2/(1 + \cos\alpha).$$

Now let $\varphi_3(\alpha)$ be chosen to satisfy

$$\Delta\varphi_3(\alpha) = \varphi_0 - \varphi_2 - \psi + C$$

$$= (1/4\pi^2)(\pi - \alpha)(1 - \cos\alpha)/\sin\alpha + (1/8\pi^2)(\pi - \alpha)^2/(1 + \cos\alpha) + C,$$

where the constant $C$ is chosen to make the right hand side have average value zero over $S^3$, and where the function $\varphi_3(\alpha)$ has a singularity at $\alpha = 0$ but not at $\alpha = \pi$.

We have yet to calculate $\varphi_3$ explicitly, but once this is done, we will have

$$\Delta_y G(V, \varphi_3)(y) = G(V, \Delta\varphi_3)(y) = G(V, \varphi_0 - \varphi_2 - \psi)(y),$$

and hence

$$\Delta_y \{A(V, \varphi_2)(y) + G(V, \varphi_3)(y)\}$$

$$= \Delta_y A(V, \varphi_2)(y) + \Delta_y G(V, \varphi_3)(y)$$

$$= \Delta_y A(V, \varphi_2)(y) + G(V, \Delta\varphi_3)(y)$$

$$= V(y) - G(V, \varphi_0 - \varphi_2 - \psi)(y) + G(V, \varphi_0 - \varphi_2 - \psi)(y)$$

$$= V(y).$$

Therefore, if we define

$$Gr(V)(y) = A(V, \varphi_2)(y) + G(V, \varphi_3)(y)$$

$$= \int_{S^3} P_{yx} V(x)\, \varphi_2(x, y)\, dx + \nabla_y \int_{S^3} P_{yx} V(x) \bullet \nabla_y \varphi_3(x, y)\, dx,$$

we have

$$\Delta_y Gr(V)(y) = V(y),$$

as desired.

This completes the proof of Theorem 3, formula (2) for the vector-valued Green's operator on $S^3$, modulo the explicit determination of the kernel function $\varphi_3$.

We break the calculation of $\varphi_3$ into three pieces.

**First piece.**

Let

$$\varphi_0(\alpha) = (-1/4\pi^2)(\pi - \alpha) \cot \alpha$$

be the fundamental solution of the Laplacian on $S^3$,

$$\Delta \varphi_0 = \delta - 1/2\pi^2.$$

The average value of $\varphi_0$ over $S^3$ is $[\varphi_0] = -1/8\pi^2$.

If we define

$$\Phi_0(\alpha) = (-1/16\pi^2) \, \alpha \, (2\pi - \alpha),$$

then we have

$$\Delta \Phi_0 = \varphi_0 - [\varphi_0].$$

**Second piece.**

Let

$$\varphi_1(\alpha) = (-1/4\pi^2)(\pi - \alpha) \csc \alpha$$

be the kernel function for the Biot-Savart operator on $S^3$ in parallel transport format. The average value of $\varphi_1$ over $S^3$ is $[\varphi_1] = -1/2\pi^2$.

If we define

$$\Phi_1(\alpha) = (-1/4\pi^2)(\pi - \alpha)(1 - \cos \alpha)/\sin \alpha,$$

then we have

$$\Delta \Phi_1 = \varphi_1 - [\varphi_1].$$

**Third piece.**

Let

$$\psi(\alpha) = (1/8\pi^2)(\pi - \alpha)^2 / (1 + \cos\alpha),$$

which is the right-most term in the expression for $\Delta\varphi_3$ in the statement of Theorem 3, formula (2).

The average value of $\psi$ over $S^3$ is

$$[\psi] = -1/2\pi^2 + 1/12.$$

If we find the function $\Psi(\alpha)$ which satisfies

$$\Delta\Psi = \psi - [\psi],$$

and which has no singularity at $\alpha = \pi$, then

$$\varphi_3 = \Phi_0 - \Phi_1 + \Psi$$

satisfies

$$\Delta\varphi_3 = (\varphi_0 - [\varphi_0]) - (\varphi_1 - [\varphi_1]) + (\psi - [\psi]),$$

and is the kernel function we are seeking.

## 31. Finding an explicit formula for $\Psi$.

The core issue is to solve the equation

$$\Delta F(\alpha) = (\pi - \alpha)^2 / (1 + \cos \alpha),$$

after which we will have

$\Psi = (1/8\pi^2) F +$ an appropriate function whose Laplacian is constant.

Note that

$$\Delta F = F'' + 2 F' \cot \alpha = (1/\sin^2\alpha) (\sin^2\alpha \, F')',$$

which tells us that $\sin^2\alpha$ is an integrating factor for our equation.

Thus we must solve

$$(\sin^2\alpha \, F')' = (\pi - \alpha)^2 \sin^2\alpha / (1 + \cos \alpha) = (\pi - \alpha)^2 (1 - \cos \alpha).$$

This integrates to

$$\sin^2\alpha \, F' = -(\pi - \alpha)^3/3 - (\pi - \alpha)^2 \sin \alpha + 2(\pi - \alpha) \cos \alpha + 2 \sin \alpha.$$

Notice that the right hand side vanishes at $\alpha = \pi$, so we do not add a constant of integration. Dividing through by $\sin^2\alpha$, we get

$$F' = -(\pi - \alpha)^3 / 3 \sin^2\alpha - (\pi - \alpha)^2 / \sin \alpha + \big(2(\pi - \alpha) \cos \alpha + 2 \sin \alpha\big) / \sin^2\alpha.$$

The right-hand half of this expression,

$$(2(\pi - \alpha) \cos \alpha + 2 \sin \alpha) / \sin^2 \alpha ,$$

blows up as $\alpha \to 0$, but approaches $0$ as $\alpha \to \pi$, and integrates to

$$-2 (\pi - \alpha) / \sin \alpha .$$

The left-hand half of this expression,

$$\omega(\alpha) = -(\pi - \alpha)^3 / 3 \sin^2 \alpha - (\pi - \alpha)^2 / \sin \alpha ,$$

blows up as $\alpha \to 0$, but approaches $0$ as $\alpha \to \pi$.

As mentioned earlier, the integral of $\omega(\alpha)$ is a non-elementary function involving polylogarithms, so for momentary convenience we define

$$\Omega(\alpha) = - \int_\alpha^\pi \omega(\alpha) \, d\alpha .$$

Thus $\Omega'(\alpha) = \omega(\alpha)$ and $\Omega(\pi) = 0$.

Then our equation for $F'$,

$$F' = -(\pi - \alpha)^3 / 3 \sin^2 \alpha - (\pi - \alpha)^2 / \sin \alpha + (2(\pi - \alpha) \cos \alpha + 2 \sin \alpha) / \sin^2 \alpha .$$

$$= \omega(\alpha) + (2(\pi - \alpha) \cos \alpha + 2 \sin \alpha) / \sin^2 \alpha ,$$

integrates to

$$F(\alpha) = \Omega(\alpha) - 2 (\pi - \alpha) / \sin \alpha .$$

This function $F(\alpha)$ satisfies the equation

$$\Delta F(\alpha) = (\pi - \alpha)^2 / (1 + \cos \alpha) .$$

Since the function $G(\alpha) = \frac{1}{2}(\pi - \alpha)\cot\alpha$ has Laplacian 1 and no singularity at $\alpha = \pi$, we conclude that the function

$$\Psi(\alpha) = (1/8\pi^2)\, F(\alpha) - (-1/2\pi^2 + 1/12)\, G(\alpha)$$

$$= (1/8\pi^2)\, (\Omega(\alpha) - 2(\pi - \alpha)/\sin\alpha) + (1/4\pi^2 - 1/24)(\pi - \alpha)\cot\alpha$$

satisfies the equation

$$\Delta\Psi = \psi - [\psi],$$

and has no singularity at $\alpha = \pi$, as desired.

## 32. Explicit formula for $\varphi_3$.

Putting this all together, we have that

$$\varphi_3 = \Phi_0 - \Phi_1 + \Psi$$

$$= (-1/16\pi^2)\,\alpha\,(2\pi - \alpha) + (1/4\pi^2)(\pi - \alpha)(1 - \cos\alpha)/\sin\alpha$$

$$+ (1/8\pi^2)(\Omega(\alpha) - 2(\pi - \alpha)/\sin\alpha) + (1/4\pi^2 - 1/24)(\pi - \alpha)\cot\alpha$$

$$= (-1/24)(\pi - \alpha)\cot\alpha - (1/16\pi^2)\,\alpha\,(2\pi - \alpha) + (1/8\pi^2)\,\Omega(\alpha),$$

$$= (-1/24)(\pi - \alpha)\cot\alpha - (1/16\pi^2)\,\alpha\,(2\pi - \alpha)$$

$$+ (1/8\pi^2) \int_\alpha^\pi \left((\pi - \alpha)^3/(3\sin^2\alpha)\right) + \left((\pi - \alpha)^2/\sin\alpha\right)\, d\alpha,$$

thanks to two cancellations.

This completes the proof of Theorem 3, formula (2).

# D. PROOFS IN $H^3$ IN PARALLEL TRANSPORT FORMAT

## 33. The hyperboloid model of hyperbolic 3-space.

In part C, we considered the standard Euclidean inner product in $R^4$,

$$<x,y> = x_0 y_0 + x_1 y_1 + x_2 y_2 + x_3 y_3,$$

and then focused on the unit 3-sphere

$$S^3 = \{x \in R^4 : <x,x> = 1\}.$$

This concrete model for the 3-sphere of constant curvature $+1$ is so common that we used it without further comment.

We now replace $R^4$ by Minkowski space $R^{1,3}$ with the indefinite inner product

$$<x,y> = x_0 y_0 - x_1 y_1 - x_2 y_2 - x_3 y_3,$$

and then regard hyperbolic 3-space as the hyperboloid

$$H^3 = \{x \in R^{1,3} : <x,x> = 1 \text{ and } x_0 > 0\}.$$

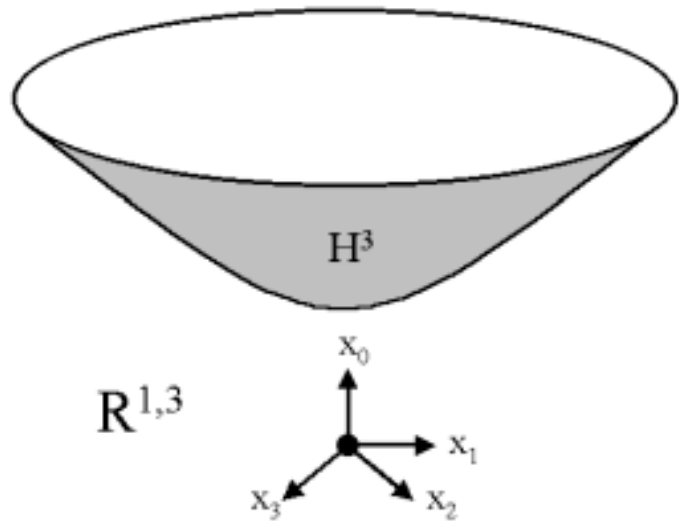

**The hyperboloid model of hyperbolic 3-space
inside Minkowski 4-space**

This concrete model for the hyperbolic 3-space of constant curvature $-1$ is not nearly so common as the Poincaré ball model or the upper half-space model. But the hyperboloid model lets us exploit the interplay between the geometry of $H^3$ and the linear algebra of $R^{1,3}$ in much the same way that we did in part C for $S^3$ and $R^4$, and hence lets us give proofs of the formulas for the Biot-Savart and Green's operators in $H^3$ which are very closely patterned on the corresponding proofs already given for $S^3$.

Since the inner product $<x,x>$ is positive for all $x \in H^3$, these radius vectors to $H^3$ are "timelike". Moreover, if $x, y \in H^3$, then we have $<x,y> \geq 1$, with equality if and only if $x = y$.

We will freely identify points in $R^{1,3}$ with the vectors at the origin that point to them, and we will also identify the tangent spaces at various points in $R^{1,3}$ in the usual way, so that we can say $v \in T_x H^3$ if and only if $<x,v> = 0$. We necessarily have $<v,v> < 0$ for all nonzero $v \in T_x H^3$; in other words, the tangent vectors to $H^3$ are "spacelike".

We will always write $<v,w>$ for the inner product of vectors when we're thinking of them as belonging to $R^{1,3}$. If $v$ and $w$ both lie in $T_x H^3$, we define

$$v \bullet w = -<v,w>.$$

Because the tangent vectors to $H^3$ are spacelike, this gives a positive definite inner product on each tangent space $T_x H^3$. The resulting Riemannian metric is complete and has constant curvature $-1$.

## 34. Isometries of $H^3$.

A linear transformation $A: R^{1,3} \to R^{1,3}$ is called an *isometry* if

$$< A(x), A(y) > \; = \; < x, y >$$

for all $x$ and $y \in R^{1,3}$.

An isometry of $R^{1,3}$ takes the two-sheeted hyperboloid $<x, x> = 1$ to itself, but may interchange the two sheets. If it does not interchange these sheets, then in particular it takes the "upper" sheet $H^3$ to itself. In that case, it also preserves the Riemannian metric on $H^3$, so is an isometry of $H^3$ in the sense of Riemannian geometry.

One example of an isometry of $R^{1,3}$ which preserves $H^3$ is given by the linear transformation $A$ with matrix

$$\begin{matrix} \cosh t & \sinh t & 0 & 0 \\ \sinh t & \cosh t & 0 & 0 \\ 0 & 0 & 1 & 0 \\ 0 & 0 & 0 & 1 \end{matrix}.$$

A family of examples of isometries of $R^{1,3}$ which preserve $H^3$ consists of those linear transformations $A$ which fix identically the $x_0$-axis, and act as Euclidean rigid motions on the $x_1 x_2 x_3$-space.

Compositions of these examples provide all possible isometries of $R^{1,3}$ which take $H^3$ to itself, and in particular provide all possible isometries of $H^3$ in the sense of Riemannian geometry.

Half of these are orientation preserving, and constitute the group $SO(1,3)$.

# 35. Geodesics in $H^3$.

Let $x \in H^3$, so that $<x,x> = 1$, and suppose that $v$ is a unit vector tangent to $H^3$ at $x$, so that $<x,v> = 0$ and $<v,v> = -1$. Consider the curve

$$G(t) = \cosh t \ x + \sinh t \ v,$$

which we easily check to lie in $H^3$ and to be parametrized by arc length.

If we complete $v$ to an orthonormal basis $v, u, w$ of $T_x H^3$, then the linear map of $R^{1,3}$ that fixes $x$ and $v$ but sends $u$ and $w$ to $-u$ and $-w$ is an isometry of $H^3$ whose fixed point set is the curve $G$. As usual in Riemannian geometry, this tells us that $G$ is a geodesic.

The inner product $<x, G(t)> = \cosh t$ is the hyperbolic cosine of the distance $t$ in $H^3$ from $x$ to $G(t)$. More generally, if $x$ and $y$ are any two points of $H^3$, the inner product $<x, y>$ is the hyperbolic cosine of the distance in $H^3$ between them.

Now let $x$ and $y$ be any pair of distinct points in $H^3$. Then the vector $v = y - <x, y> x$ is nonzero and orthogonal to $x$, so that $v \in T_x H^3$. Moreover, $|v|^2 = <v, v> = 1 - <x, y>^2 < 0$, so $v \bullet v = <x, y>^2 - 1 > 0$, and hence $v / |v|$ is a unit vector tangent to $H^3$ at $x$. Therefore the unique, length-minimizing geodesic from $x$ to $y$ is given by

$$G(t) = \cosh t \ x + \sinh t \ (v / |v|).$$

If $\alpha$ is the distance from $x$ to $y$ in $H^3$, so that $\cosh \alpha = <x, y>$, then we can rewrite this as

$$G(t) = \cosh t \ x + \sinh t \ ((y - \cosh \alpha \ x) / \sinh \alpha).$$

When $t = \alpha$, we have $G(\alpha) = y$ and

$$G'(\alpha) = \sinh \alpha \ x + \cosh \alpha \ ((y - \cosh \alpha \ x) / \sinh \alpha)$$

$$= (\cosh \alpha \ y - x) / \sinh \alpha.$$

Since $G$ is parametrized by arc-length and is length-minimizing, and $\alpha$ is the distance in $H^3$ from $x$ to $y$, we deduce that

$$\nabla_y \alpha(x, y) = G'(\alpha) = (\cosh \alpha \ y - x) / \sinh \alpha.$$

We note that
$$P_{xy} \nabla_y \alpha(x, y) = -\nabla_x \alpha(x, y).$$

# 36. The vector triple product in $R^{1,3}$.

Let $X_0 = (1, 0, 0, 0)$, $X_1 = (0, 1, 0, 0)$, $X_2 = (0, 0, 1, 0)$ and $X_3 = (0, 0, 0, 1)$ be the standard orthonormal basis for $R^{1,3}$. Then for $x, v, w \in R^{1,3}$, we define the *vector triple product*

$$[x, v, w] = \det \begin{pmatrix} X_0 & X_1 & X_2 & X_3 \\ v_0 & v_1 & v_2 & v_3 \\ w_0 & w_1 & w_2 & w_3 \\ X_0 & -X_1 & -X_2 & -X_3 \end{pmatrix} \in R^{1,3}.$$

Note that $[x, v, w]$ is an alternating multilinear function of its three arguments.

For $y \in R^{1,3}$, we have

$$<[x, v, w], y> = \det \begin{pmatrix} x_0 & x_1 & x_2 & x_3 \\ v_0 & v_1 & v_2 & v_3 \\ w_0 & w_1 & w_2 & w_3 \\ y_0 & y_1 & y_2 & y_3 \end{pmatrix}.$$

Thus $[x, v, w]$ is a vector in $R^{1,3}$ whose inner product with each of $x$, $v$ and $w$ is zero.

If $A \in SO(1,3)$, we have

$$[Ax, Av, Aw] = A[x, v, w].$$

## 37. Explicit formula for parallel transport in $H^3$.

Let $x$ and $y$ be any pair of distinct points in $H^3 \subset R^{1,3}$. Then there is a unique element $M \in SO(1,3)$ that maps $x$ to $y$ and leaves fixed the two-dimensional subspace of $R^{1,3}$ consisting of all vectors orthogonal to both $x$ and $y$. For $v \in R^{1,3}$, this mapping is given by

$$M(v) = v - \frac{<x+y, v>}{1 + <x, y>} x + \left\langle \frac{1 + 2<x, y>}{1 + <x, y>} x - \frac{1}{1 + <x, y>} y, v \right\rangle y.$$

Direct calculation shows that $M \in SO(1,3)$, and is therefore the desired isometry.

If $v \in T_x H^3$, then $<v, x> = 0$, and the formula for $M(v)$ simplifies to

$$M(v) = v - \frac{<y, v>}{1 + <x, y>} (x + y).$$

In this case, $M(v)$ is the result of parallel transport of $v$ from $x$ to $y$, so we will write

$$P_{yx} v = M(v) = v - \frac{<y, v>}{1 + <x, y>} (x + y).$$

Note that the formulas in this section are the same as those in section 16 in the case of $S^3 \subset R^4$, except for the change in the definition of the inner product.

# 38. The method of moving frames in $H^3$.

The computations which lead to the various integral formulas in $H^3$ in parallel transport format will be carried out using "moving frames" which are obtained via parallel transport from a single orthonormal frame, just as they were on $S^3$, as follows.

Let $x = (1, 0, 0, 0)$ in $H^3 \subset R^{1,3}$, and let $X_1$, $X_2$, $X_3$ be the usual orthonormal basis for the tangent space $T_xH^3$ to $H^3$ at $x$.

For each point $y = (y_0, y_1, y_2, y_3) \in H^3$, we parallel transport this orthonormal basis for $T_xH^3$ along the unique geodesic from $x$ to $y$ to obtain an orthonormal basis

$$E_1 = P_{yx}X_1, \quad E_2 = P_{yx}X_2, \quad E_3 = P_{yx}X_3$$

for $T_yH^3$.

Using the formula for $P_{yx}V$, we calculate that

$$E_1 = (y_1, \ 1 + y_1^2/(1 + y_0), \ y_1y_2/(1 + y_0), \ y_1y_3/(1 + y_0))$$

$$E_2 = (y_2, \ y_1y_2/(1 + y_0), \ 1 + y_2^2/(1 + y_0), \ y_2y_3/(1 + y_0))$$

$$E_3 = (y_3, \ y_1y_3/(1 + y_0), \ y_2y_3/(1 + y_0), \ 1 + y_3^2/(1 + y_0)).$$

Next, let $\theta_1$, $\theta_2$, $\theta_3$ denote the 1-forms on $H^3$ which are dual to $E_1$, $E_2$, $E_3$, in the sense that $\theta_i(V) = E_i \bullet V = -<E_i, V>$ for all vectors $V$ and for $i = 1, 2, 3$.

We can write

$$\theta_1 = -y_1 \, dy_0 + (1 + y_1^2/(1 + y_0)) \, dy_1 + y_1y_2/(1 + y_0) \, dy_2 + y_1y_3/(1 + y_0) \, dy_3$$

$$\theta_2 = -y_2 \, dy_0 + y_1y_2/(1 + y_0) \, dy_1 + (1 + y_2^2/(1 + y_0)) \, dy_2 + y_2y_3/(1 + y_0) \, dy_3$$

$$\theta_3 = -y_3 \, dy_0 + y_1y_3/(1 + y_0) \, dy_1 + y_2y_3/(1 + y_0) \, dy_2 + (1 + y_3^2/(1 + y_0)) \, dy_3,$$

where we understand the right hand sides as restricted to $H^3$.

Since $y_0^2 - y_1^2 - y_2^2 - y_3^2 = 1$ on $H^3$, we can add a fourth equation to the three above,

$$0 = y_0 \, dy_0 - y_1 \, dy_1 - y_2 \, dy_2 - y_3 \, dy_3 \, .$$

Using this fourth equation, we can simplify the first three to read

$$\theta_1 = -y_1/(1 + y_0) \, dy_0 + dy_1$$

$$\theta_2 = -y_2/(1 + y_0) \, dy_0 + dy_2$$

$$\theta_3 = -y_3/(1 + y_0) \, dy_0 + dy_3 \, .$$

We then invert this system to obtain

$$dy_0 = y_1 \, \theta_1 + y_2 \, \theta_2 + y_3 \, \theta_3$$

$$dy_1 = \theta_1 + (y_1/(1 + y_0)) \, (y_1 \, \theta_1 + y_2 \, \theta_2 + y_3 \, \theta_3)$$

$$dy_2 = \theta_2 + (y_2/(1 + y_0)) \, (y_1 \, \theta_1 + y_2 \, \theta_2 + y_3 \, \theta_3)$$

$$dy_3 = \theta_3 + (y_3/(1 + y_0)) \, (y_1 \, \theta_1 + y_2 \, \theta_2 + y_3 \, \theta_3) \, .$$

With these equations in hand, we easily obtain

$$d\theta_1 = -y_2/(1 + y_0) \, \theta_1 \wedge \theta_2 + y_3/(1 + y_0) \, \theta_3 \wedge \theta_1$$

$$d\theta_2 = -y_3/(1 + y_0) \, \theta_2 \wedge \theta_3 + y_1/(1 + y_0) \, \theta_1 \wedge \theta_2$$

$$d\theta_3 = -y_1/(1 + y_0) \, \theta_3 \wedge \theta_1 + y_2/(1 + y_0) \, \theta_2 \wedge \theta_3 \, .$$

## 39. Curls and divergences.

If $x$ is a fixed point on $H^3$ and $V$ is a tangent vector to $H^3$ at $x$, then the parallel transport $P_{yx}V$ of $V$ to the point $y$ is defined for all $y \in H^3$, and may be viewed as a vector field on $H^3$.

We need to know the curl and divergence of this vector field.

**PROPOSITION 39.1.**
$$\nabla_y \times P_{yx}V = \frac{[y, x, V]}{1 + <x, y>} .$$

**PROPOSITION 39.2.**
$$\nabla_y \bullet P_{yx}V = \frac{-2<y, V>}{1 + <x, y>} .$$

We use the same definitions as on $S^3$ for the curl and the divergence, namely
$$\nabla \times W = \Uparrow d \downarrow W ,$$
and
$$\nabla \bullet W = * d \Downarrow W .$$

For example, suppose $W$ is the vector field $E_1 = P_{yx}X_1$. Then
$$\nabla_y \times E_1 = \Uparrow d \downarrow E_1 = \Uparrow d \theta_1$$
$$= \Uparrow - (y_2/(1 + y_0) \; \theta_1 \wedge \theta_2 + y_3/(1 + y_0) \; \theta_3 \wedge \theta_1)$$
$$= - y_2/(1 + y_0) \; E_3 + y_3/(1 + y_0) \; E_2 .$$

And likewise
$$\nabla_y \times E_2 = - y_3/(1 + y_0) \; E_1 + y_1/(1 + y_0) \; E_3$$
$$\nabla_y \times E_3 = - y_1/(1 + y_0) \; E_2 + y_2/(1 + y_0) \; E_1 .$$

Note that these formulas are the negatives of the ones obtained on $S^3$ in section 18.

They can be viewed as special instances of the formula given above in Proposition 39.1. For example, if $y = (y_0, y_1, y_2, y_3)$, $x = (1, 0, 0, 0)$ and $V = X_1$, then one easily calculates that

$$[y, x, V] = y_3 E_2 - y_2 E_3.$$

Dividing through by $1 + <x, y> = 1 + y_0$ gives the formula for $\nabla_y \times E_1$ appearing above, and so verifies Proposition 39.1 in this instance. Just as in the case of $S^3$, the verification of the proposition in this single instance also serves as its proof.

The divergence of a vector field $W$ was defined above by $\nabla \bullet W = *d \Downarrow W$. If $W$ is the vector field $E_1 = P_{yx} X_1$, then

$$\nabla \bullet E_1 = *d \Downarrow E_1 = *d(\theta_2 \wedge \theta_3)$$

$$= *(d\theta_2 \wedge \theta_3 - \theta_2 \wedge d\theta_3)$$

$$= *(2 y_1/(1 + y_0)) \theta_1 \wedge \theta_2 \wedge \theta_3$$

$$= 2 y_1/(1 + y_0),$$

using the formulas for $d\theta_2$ and $d\theta_3$ appearing in the previous section.

The above example can be viewed as a special instance of the formula given in Proposition 39.2, since

$$<y, V> = <y, X_1> = -y_1 \text{ and } 1 + <x, y> = 1 + y_0.$$

And just as before, this single verification of the proposition also serves as its proof.

## 40. Statement of the Key Lemma.

The following result plays a lead role in obtaining the formulas for the Biot-Savart and Green's operators on $H^3$ in parallel transport format.

**KEY LEMMA, hyperbolic version.**

$$\nabla_y \times \{P_{yx}V(x) \times \nabla_y \varphi\} - \nabla_y \{V(x) \bullet \nabla_x (\cosh \alpha \ \varphi)\}$$
$$= (\Delta\varphi + \varphi)(V(x) - <V(x), y> y).$$

As usual, $\alpha$ is the distance on $H^3$ between $x$ and $y$. The kernel function $\varphi(\alpha)$ is any smooth function, which typically blows up at $\alpha = 0$ and is asymptotic there to $-1/(4\pi\alpha)$.

The Key Lemma can be understood intuitively when $\varphi(\alpha)$ is the fundamental solution of the operator $\varphi \to \Delta\varphi + \varphi$ on $S^3$, in which case we will see that $\cosh \alpha \ \varphi(\alpha)$ is, up to an additive constant, the fundamental solution of the scalar Laplacian there.

Then, just as on $S^3$, the Key Lemma can be viewed as a pointwise or infinitesimal version of Maxwell's equation

$$\nabla \times B - \partial E/\partial t \ = V \ .$$

## 41. The first term on the LHS of the Key Lemma.

The plan is the same as it was for $S^3$. Applying an appropriate isometry of $H^3$, we can assume that $x$ is the point $(1, 0, 0, 0)$ and that, after scaling, $V(x) = X_1 = (0, 1, 0, 0)$. At each point $y = (y_0, y_1, y_2, y_3)$ of $H^3$, we will expand both sides of the Key Lemma in terms of the orthonormal basis $E_1$, $E_2$ and $E_3$ for $T_y H^3$, and prove equality component-wise.

We begin with the first term on the LHS of the equation, and note that

$$P_{yx} V(x) \times \nabla_y \varphi = P_{yx} V(x) \times (-P_{yx} \nabla_x \varphi) = P_{yx} (- V(x) \times \nabla_x \varphi).$$

From section 35, we have that

$$\nabla_x \alpha = (\cosh \alpha \; x - y) / \sinh \alpha$$

$$= ((y_0, 0, 0, 0) - (y_0, y_1, y_2, y_3)) / \sinh \alpha$$

$$= (-1 / \sinh \alpha)(y_1 X_1 + y_2 X_2 + y_3 X_3).$$

Hence

$$\nabla_x \varphi = \varphi'(\alpha) \nabla_x \alpha = (- \varphi'(\alpha) / \sinh \alpha)(y_1 X_1 + y_2 X_2 + y_3 X_3).$$

Thus

$$- V(x) \times \nabla_x \varphi = (-X_1) \times (- \varphi'(\alpha) / \sinh \alpha)(y_1 X_1 + y_2 X_2 + y_3 X_3)$$

$$= (\varphi'(\alpha) / \sinh \alpha)(-y_3 X_2 + y_2 X_3).$$

Hence

$$P_{yx} V(x) \times \nabla_y \varphi = P_{yx} (- V(x) \times \nabla_x \varphi)$$

$$= (\varphi'(\alpha) / \sinh \alpha)(-y_3 E_2 + y_2 E_3).$$

Now we take curls and evaluate the first term on the LHS of the Key Lemma,

$$\nabla_y \times \{P_{yx}V(x) \times \nabla_y \varphi\} = \nabla_y \times \{(\varphi'(\alpha) / \sinh \alpha) (-y_3 E_2 + y_2 E_3)\} .$$

We begin by computing the curl of the vector field

$$W = W(y) = -y_3 E_2 + y_2 E_3 ,$$

following the recipe $\nabla_y \times W = \Uparrow d \Downarrow W$. The computation is similar to the one on $S^3$, and yields

$$\nabla_y \times W = 2 (y_0 - y_1^2/(1 + y_0)) E_1 - 2 (y_1 y_2 / (1 + y_0)) E_2 - 2 (y_1 y_3 / (1 + y_0)) E_3 .$$

Now to compute $\nabla_y \times \{(\varphi'(\alpha) / \sinh \alpha) (-y_3 E_2 + y_2 E_3)\}$, we use the formula

$$\nabla \times (f W) = f (\nabla \times W) + \nabla f \times W ,$$

with $f(\alpha) = \varphi'(\alpha) / \sinh \alpha$, and find that

$$\nabla_y \times \{P_{yx}V(x) \times \nabla_y \varphi\} = \nabla_y \times \{(\varphi'(\alpha) / \sinh \alpha) (-y_3 E_2 + y_2 E_3)\}$$

$$= \{ (2\varphi'/\sinh \alpha) (y_0 - y_1^2/(1 + y_0)) + (\varphi''/\sinh^2\alpha - \varphi' \cosh \alpha/\sinh^3\alpha) (y_2^2 + y_3^2) \} E_1$$

$$+ \{ (-2\varphi'/\sinh \alpha) y_1 y_2/(1 + y_0) + (\varphi''/\sinh^2\alpha - \varphi' \cosh \alpha/\sinh^3\alpha) (-y_1 y_2) \} E_2$$

$$+ \{ (-2\varphi'/\sinh \alpha) y_1 y_3/(1 + y_0) + (\varphi''/\sinh^2\alpha - \varphi' \cosh \alpha/\sinh^3\alpha) (-y_1 y_3) \} E_3 .$$

This puts the first term on the LHS of the Key Lemma in the form we want it.

## 42. The second term on the LHS of the Key Lemma.

Now, just as we did for the first term, we want to express the second term

$$- \nabla_y \{V(x) \bullet \nabla_x (\cosh \alpha \; \varphi)\}$$

on the LHS of the Key Lemma as a linear combination of the orthonormal basis vectors $E_1$, $E_2$ and $E_3$ for $T_y H^3$.

We follow the procedure used on $S^3$ in section 21, and find that

$$- \nabla_y \{V(x) \bullet \nabla_x (\cosh \alpha \; \varphi)\}$$

$$= \{ (\varphi'' \cosh \alpha / \sinh^2 \alpha + \varphi'(1/\sinh \alpha - 1/\sinh^3 \alpha)) y_1^2$$

$$+ (\varphi' \cosh \alpha / \sinh \alpha + \varphi) (1 + y_1^2/(1 + y_0)) \} E_1$$

$$+ \{ (\varphi'' \cosh \alpha / \sinh^2 \alpha + \varphi'(1/\sinh \alpha - 1/\sinh^3 \alpha)) y_1 y_2$$

$$+ (\varphi' \cosh \alpha / \sinh \alpha + \varphi) y_1 y_2/(1 + y_0) \} E_2$$

$$+ \{ (\varphi'' \cosh \alpha / \sinh^2 \alpha + \varphi'(1/\sinh \alpha - 1/\sinh^3 \alpha)) y_1 y_3$$

$$+ (\varphi' \cosh \alpha / \sinh \alpha + \varphi) y_1 y_3/(1 + y_0) \} E_3 \, .$$

If we combine the two terms on the LHS of the Key Lemma and simplify, we get

$$\nabla_y \times \{P_{yx} V(x) \times \nabla_y \varphi\} - \nabla_y \{V(x) \bullet \nabla_x (\cosh \alpha \; \varphi)\}$$

$$= \{(\varphi'' + 2 \varphi' \cosh \alpha / \sinh \alpha + \varphi) (1 + y_1^2/(1 + y_0)) \} E_1$$

$$+ \{(\varphi'' + 2 \varphi' \cosh \alpha / \sinh \alpha + \varphi) (y_1 y_2/(1 + y_0)) \} E_2$$

$$+ \{(\varphi'' + 2 \varphi' \cosh \alpha / \sinh \alpha + \varphi) (y_1 y_3/(1 + y_0)) \} E_3 \, .$$

## 43. The RHS of the Key Lemma.

Now we express the RHS

$$(\Delta\varphi + \varphi)(V(x) - <V(x), y> y)$$

as a linear combination of $E_1$, $E_2$ and $E_3$.

Recall that

$$\Delta\varphi = \varphi'' + 2\varphi' \cosh\alpha/\sinh\alpha,$$

and hence

$$\Delta\varphi + \varphi = \varphi'' + 2\varphi' \cosh\alpha/\sinh\alpha + \varphi.$$

We follow the procedure used on $S^3$ in section 22, and find that

$$(\Delta\varphi + \varphi)(V(x) - <V(x), y> y)$$

$$= \{(\varphi'' + 2\varphi' \cosh\alpha/\sinh\alpha + \varphi)(1 + y_1^2/(1+y_0))\} E_1$$

$$+ \{(\varphi'' + 2\varphi' \cosh\alpha/\sinh\alpha + \varphi)(y_1 y_2/(1+y_0))\} E_2$$

$$+ \{(\varphi'' + 2\varphi' \cosh\alpha/\sinh\alpha + \varphi)(y_1 y_3/(1+y_0))\} E_3.$$

This matches the LHS, as given in the previous section, and so completes the proof of the Key Lemma.

## 44. Maxwell's equation in $H^3$.

We now apply the Key Lemma to the proof of Maxwell's equation

$$\nabla \times B = V + \partial E/\partial t$$

in $H^3$. In the Key Lemma, we use the kernel function

$$\varphi(\alpha) = (-1/4\pi) \operatorname{csch} \alpha,$$

which, as we noted earlier, is a fundamental solution of the shifted Laplacian on $H^3$,

$$\Delta\varphi + \varphi = \delta.$$

With this choice of kernel function, we have that

$$(\cosh \alpha)\, \varphi(\alpha) = (-1/4\pi) \coth \alpha,$$

which differs from the fundamental solution of the Laplacian on $H^3$,

$$\varphi_0(\alpha) = (-1/4\pi)(\coth \alpha - 1),$$

only by a constant. Since we take the gradient of this function below, the difference is irrelevant.

Hence with these choices, the Key Lemma reads

$$\nabla_y \times \{P_{yx} V(x) \times \nabla_y \varphi\} - \nabla_y \{V(x) \bullet \nabla_x \varphi_0\}$$
$$= \delta(x, y)\, (V(x) - <V(x), y>y).$$

Keeping these choices of kernel functions, we now define the Biot-Savart operator by formula (3) of Theorem 2,

$$BS(V)(y) = \int_{H^3} P_{yx} V(x) \times \nabla_y \varphi(x, y)\, dx,$$

for smooth, compactly supported vector fields $V$.

Then we integrate the Key Lemma over $H^3$ and get

$$\int_{H^3} \nabla_y \times \{P_{yx}V(x) \times \nabla_y \varphi\} \, dx \; - \; \int_{H^3} \nabla_y \{V(x) \bullet \nabla_x \varphi_0\} \, dx$$

$$= \int_{H^3} \delta(x, y) \left(V(x) - <V(x), y> y\right) dx \,.$$

We process this just as we did on $S^3$ in section 23, and obtain

$$\nabla_y \times BS(V)(y) \; = \; V(y) \; + \; \nabla_y \int_{H^3} V(x) \bullet \nabla_x \varphi_0(x, y) \, dx \,.$$

As in $R^3$ and on $S^3$, we recognize this as Maxwell's equation

$$\nabla \times B \; = \; V \; + \; \partial E/\partial t \,.$$

## 45. Proof of Theorem 2, formula (3).

Let $V$ be a smooth compactly supported vector field in $H^3$. Thinking of $V$ as a steady current flow, we defined the corresponding magnetic field $BS(V)$ by

$$BS(V)(y) \; = \; \int_{H^3} P_{yx}V(x) \times \nabla_y \varphi(x, y) \, dx \,,$$

where $\varphi(\alpha) \; = \; (-1/4\pi) \, \text{csch} \, \alpha$.

We saw in the preceding section that, as a consequence of the Key Lemma, the field $BS(V)$ satisfies Maxwell's equation.

To complete the proof of Theorem 2, formula (3), we need only show that $BS(V)$ is divergence-free.

But the argument for this in $H^3$ is identical to that in $S^3$, so we are done.

# 46. The Green's operator on $H^3$ in parallel transport format.

We turn now to the derivation of formula (3) of Theorem 3 for the vector-valued Green's operator on $H^3$,

$$Gr(V)(y) = \int_{H^3} P_{yx}V(x)\, \varphi_2(x, y)\, dx + \nabla_y \int_{H^3} P_{yx}V(x) \bullet \nabla_y \varphi_3(x, y)\, dx,$$

where

$$\varphi_2(\alpha) = (-1/4\pi)\operatorname{csch}\alpha + (1/4\pi)\,\alpha/(1 + \cosh\alpha)$$

and

$$\varphi_3(\alpha) = (1/4\pi)\,\alpha/(e^{2\alpha} - 1) + (1/4\pi) \int_0^\alpha \left((\alpha/\sinh\alpha) - (\alpha^2/2\sinh^2\alpha)\right) d\alpha.$$

The naming of the kernel functions $\varphi_2$ and $\varphi_3$ above leaves room in the $\varphi$-family for the fundamental solution of the Laplacian on $H^3$,

$$\varphi_0(\alpha) = (-1/4\pi)(\coth\alpha - 1) \qquad\qquad \Delta\varphi_0 = \delta$$

and for the kernel function of the Biot-Savart operator,

$$\varphi_1(\alpha) = (-1/4\pi)\operatorname{csch}\alpha \qquad\qquad \Delta\varphi_1 + \varphi_1 = \delta.$$

## 47. Plan of the proof.

Consider the usual vector convolutions,

$$A(V, \varphi)(y) = \int_{H^3} P_{yx} V(x) \, \varphi(x, y) \, dx$$

$$B(V, \varphi)(y) = \int_{H^3} P_{yx} V(x) \times \nabla_y \varphi(x, y) \, dx$$

$$g(V, \varphi)(y) = \int_{H^3} P_{yx} V(x) \bullet \nabla_y \varphi(x, y) \, dx$$

$$G(V, \varphi)(y) = \nabla_y \, g(V, \varphi)(y) = \nabla_y \int_{H^3} P_{yx} V(x) \bullet \nabla_y \varphi(x, y) \, dx \, .$$

***Step 1.*** First we will determine the unknown kernel function $\varphi$ so as to make

$$\nabla_y \times A(V, \varphi)(y) = -B(V, \varphi_1)(y) = -BS(V)(y) \, .$$

The solution we obtain is *not* unique.

***Step 2.*** Then, as a consequence of Maxwell's equation, we will have

$$-\nabla_y \times \nabla_y \times A(V, \varphi_2)(y) = \nabla_y \times BS(V)(y)$$

$$= V(y) - \nabla_y \int_{H^3} \nabla_x \bullet V(x) \, \varphi_0(x, y) \, dx \, .$$

Hence

$$\Delta_y \, A(V, \varphi_2)(y) = -\nabla_y \times \nabla_y \times A(V, \varphi_2)(y) + \nabla_y(\nabla_y \bullet A(V, \varphi_2)(y))$$

$$= V(y) - \nabla_y \int_{H^3} \nabla_x \bullet V(x) \, \varphi_0(x, y) \, dx \, + \, \nabla_y(\nabla_y \bullet A(V, \varphi_2)(y))$$

$$= V(y) + \text{some gradient field} \, .$$

We then find the kernel function $\varphi_3$ so that $\Delta_y \, G(V, \varphi_3)(y) = G(V, \Delta\varphi_3)(y)$ is the negative of the above gradient field, in which case, defining

$$Gr(V)(y) = A(V, \varphi_2)(y) + G(V, \varphi_3)(y)$$

yields the desired equation $\Delta_y \, Gr(V)(y) = V(y) \, .$

## 48. Some formulas that we will use.

In the formulas below, $V$ is a single tangent vector at the single point $x \in H^3$.

(1) $$P_{yx}V = V - \frac{\langle y, V \rangle}{1 + \langle x, y \rangle}(x + y)$$

(2) $$\nabla_y \times P_{yx}V = \frac{[y, x, V]}{1 + \langle x, y \rangle}$$

(3) $$\nabla_y \bullet P_{yx}V = \frac{-2\langle y, V \rangle}{1 + \langle x, y \rangle}$$

(4) $$P_{yx}V \times \nabla_y \varphi = -[y, x, V]\,\varphi'(\alpha) / \sinh \alpha \ .$$

Formula (1) appears in section 37, while formulas (2) and (3) appear in section 39.

Formula (4) in the special case that $x = (1, 0, 0, 0)$ and $V = X_1 = (0, 1, 0, 0)$ is obtained from the formula

$$P_{yx}V \times \nabla_y \varphi = (\varphi'(\alpha) / \sinh \alpha)(-y_3 E_2 + y_2 E_3)$$

of section 41, together with the remark of section 39 that in this case we have

$$[y, x, V] = y_3 E_2 - y_2 E_3 \ .$$

We leave it to the reader to confirm, as we did earlier in similar circumstances, that this single case of formula (4) can serve as its proof.

We note that formulas (1), (2) and (3) are identical to the corresponding formulas for $S^3$ in section 27, but that formula (4) has the opposite sign, in addition to the expected replacement of "$\sin \alpha$" by "$\sinh \alpha$".

## 49. The curl of $A(V, \varphi)(y)$.

$\nabla_y \times A(V, \varphi)(y) = \int_{H^3} \nabla_y \times \{P_{yx}V(x) \, \varphi(x, y)\} \, dx$

$= \int_{H^3} \{\nabla_y \times P_{yx}V(x)\} \, \varphi(x, y) \; - \; P_{yx}V(x) \times \nabla_y \varphi(x, y) \, dx$

$= \int_{H^3} \dfrac{[y, x, V(x)]}{1 + <x, y>} \varphi(x, y) \; + \; [y, x, V(x)] \dfrac{\varphi'(\alpha)}{\sinh \alpha} \, dx$

$= \int_{H^3} [y, x, V(x)] \left\{ \dfrac{\varphi}{1 + \cosh \alpha} + \dfrac{\varphi'}{\sinh \alpha} \right\} dx \; ,$

thanks to formulas (2) and (4) above.

On the other hand,

$BS(V)(y) = B(V, \varphi_1)(y) = \int_{H^3} P_{yx}V(x) \times \nabla_y \varphi_1(x, y) \, dx$

$= \int_{H^3} [y, x, V(x)] \dfrac{-\varphi_1'(\alpha)}{\sinh \alpha} \, dx \; ,$

again by formula (4).

Thus, in order that

$$\nabla_y \times A(V, \varphi)(y) = -BS(V)(y) = -B(V, \varphi_1)(y) \; ,$$

we must choose $\varphi$ to satisfy the ODE

$$\dfrac{\varphi}{1 + \cosh \alpha} + \dfrac{\varphi'}{\sinh \alpha} = \dfrac{\varphi_1'}{\sinh \alpha} \; .$$

This ODE is solved just like the corresponding one in section 28, and we get

$$(1 + \cosh \alpha) \, \varphi = (1 + \cosh \alpha) \, \varphi_1 + (1/4\pi) (\alpha + C) \; .$$

The choice of the constant $C$ represents the nonuniqueness in carrying out Step 1. We discuss this further in section 52, but here set $C = 0$ and continue.

Dividing through by $(1 + \cosh \alpha)$, we get

$$\varphi(\alpha) = \varphi_1(\alpha) + (1/4\pi) \alpha / (1 + \cosh \alpha)$$

$$= (-1/4\pi) \operatorname{csch} \alpha + (1/4\pi) \alpha / (1 + \cosh \alpha).$$

Note that $\varphi(\alpha)$ has a singularity at $\alpha = 0$.

We call this solution $\varphi_2(\alpha)$, and with it we have completed Step 1:

$$\nabla_y \times A(V, \varphi_2)(y) = -B(V, \varphi_1)(y) = -BS(V)(y).$$

We take the negative curl of both sides to get

$$-\nabla_y \times \nabla_y \times A(V, \varphi_2)(y) = \nabla_y \times BS(V)(y)$$

$$= V(y) - \nabla_y \int_{H^3} (\nabla_x \bullet V(x)) \varphi_0(x, y) \, dx,$$

thanks to Maxwell's equation. Integrating by parts, we have

$$\int_{H^3} (\nabla_x \bullet V(x)) \varphi_0(x, y) \, dx = -\int_{H^3} V(x) \bullet \nabla_x \varphi_0(x, y) \, dx$$

$$= \int_{H^3} P_{yx} V(x) \bullet \nabla_y \varphi_0(x, y) \, dx$$

$$= g(V, \varphi_0)(y).$$

Hence

$$\Delta_y A(V, \varphi_2)(y) = -\nabla_y \times \nabla_y \times A(V, \varphi_2)(y) + \nabla_y(\nabla_y \bullet A(V, \varphi_2)(y))$$

$$= V(y) - \nabla_y \int_{H^3} \nabla_x \bullet V(x) \varphi_0(x, y) \, dx + \nabla_y(\nabla_y \bullet A(V, \varphi_2)(y))$$

$$= V(y) - G(V, \varphi_0)(y) + \nabla_y(\nabla_y \bullet A(V, \varphi_2)(y)).$$

Just as on $S^3$, the next step is to compute $\nabla_y \bullet A(V, \varphi_2)(y)$.

## 50. The divergence of $A(V, \varphi_2)$.

We begin as follows.

$$\nabla_y \bullet A(V, \varphi_2)(y) = \int_{H^3} \nabla_y \bullet \left( P_{yx} V(x) \, \varphi_2(x, y) \right) dx$$

$$= \int_{H^3} \left( \nabla_y \bullet P_{yx} V(x) \right) \varphi_2(x, y) \, dx + \int_{H^3} P_{yx} V(x) \bullet \nabla_y \varphi_2(x, y) \, dx .$$

The second integral on the RHS above is simply $g(V, \varphi_2)(y)$, so we put that aside.

Turning to the first integral and using formula (3) and the fact that $<x, V(x)> = 0$, we follow the pattern set in section 29 for $S^3$, and get

$$\int_{H^3} \left( \nabla_y \bullet P_{yx} V(x) \right) \varphi_2(x, y) \, dx = \int_{H^3} P_{yx} V(x) \bullet \nabla_y \psi(\alpha) \, dx$$

$$= g(V, \psi)(y) ,$$

where the new kernel function $\psi$ satisfies

$$\psi'(\alpha) = [2 \sinh \alpha / (1 + \cosh \alpha)] \, \varphi_2(\alpha)$$

$$= [2 \sinh \alpha / (1 + \cosh \alpha)] \left[ (-1/4\pi) \operatorname{csch} \alpha + (1/4\pi) \alpha / (1 + \cosh \alpha) \right]$$

$$= (-1/2\pi) / (1 + \cosh \alpha) + (1/2\pi) \, \alpha \sinh \alpha / (1 + \cosh \alpha)^2 ,$$

and hence

$$\psi(\alpha) = (-1/2\pi) \, \alpha / (1 + \cosh \alpha) .$$

We have set the constant of integration to zero, so that $\psi(\alpha) \to 0$ as $\alpha \to \infty$.

Collecting information, we have

$$\nabla_y \bullet A(V, \varphi_2)(y)$$

$$= \int_{H^3} \left( \nabla_y \bullet P_{yx} V(x) \right) \varphi_2(x, y) \, dx + \int_{H^3} P_{yx} V(x) \bullet \nabla_y \varphi_2(x, y) \, dx$$

$$= g(V, \psi)(y) + g(V, \varphi_2)(y) ,$$

and hence

$$\nabla_y \left( \nabla_y \bullet A(V, \varphi_2)(y) \right) = G(V, \psi)(y) + G(V, \varphi_2)(y) .$$

## 51. Step 2 - Finding and using the kernel function $\varphi_3$.

Continuing our calculation of the Laplacian of $A(V, \varphi_2)(y)$, we have

$$\Delta_y A(V, \varphi_2)(y) = V(y) - G(V, \varphi_0)(y) + \nabla_y(\nabla_y \bullet A(V, \varphi_2)(y))$$

$$= V(y) - G(V, \varphi_0)(y) + G(V, \psi)(y) + G(V, \varphi_2)(y)$$

$$= V(y) - G(V, \varphi_0 - \varphi_2 - \psi)(y) .$$

Using the values for $\varphi_0$, $\varphi_2$ and $\psi$ given above, we calculate that

$$\varphi_0 - \varphi_2 - \psi = (1/4\pi) - (1/4\pi) \sinh \alpha / (1 + \cosh \alpha) + (1/4\pi) \alpha / (1 + \cosh \alpha) .$$

Now let $\varphi_3(\alpha)$ be chosen to satisfy

$$\Delta \varphi_3(\alpha) = \varphi_0 - \varphi_2 - \psi$$

$$= (1/4\pi) - (1/4\pi) \sinh \alpha / (1 + \cosh \alpha) + (1/4\pi) \alpha / (1 + \cosh \alpha) .$$

We have yet to calculate $\varphi_3$ explicitly, but once this is done, we will have

$$\Delta_y G(V, \varphi_3)(y) = G(V, \Delta \varphi_3)(y) = G(V, \varphi_0 - \varphi_2 - \psi)(y) ,$$

and hence

$$\Delta_y \{A(V, \varphi_2)(y) + G(V, \varphi_3)(y)\}$$

$$= \Delta_y A(V, \varphi_2)(y) + \Delta_y G(V, \varphi_3)(y)$$

$$= \Delta_y A(V, \varphi_2)(y) + G(V, \Delta \varphi_3)(y)$$

$$= V(y) - G(V, \varphi_0 - \varphi_2 - \psi)(y) + G(V, \varphi_0 - \varphi_2 - \psi)(y)$$

$$= V(y) .$$

Therefore, if we define

$$Gr(V)(y) = A(V, \varphi_2)(y) + G(V, \varphi_3)(y)$$
$$= \int_{H^3} P_{yx}V(x)\, \varphi_2(x, y)\, dx + \nabla_y \int_{H^3} P_{yx}V(x) \bullet \nabla_y \varphi_3(x, y)\, dx,$$

we have

$$\Delta_y\, Gr(V)(y) = V(y),$$

as desired.

This completes the proof of Theorem 3, formula (3) for the vector-valued Green's operator in $H^3$, modulo the explicit determination of the kernel function $\varphi_3$.

We break the calculation of $\varphi_3$ into two pieces, and begin by writing

$$\varphi_0 - \varphi_2 - \psi = (1/4\pi) - (1/4\pi) \sinh \alpha / (1 + \cosh \alpha) + (1/4\pi) \alpha / (1 + \cosh \alpha).$$

$$= \qquad \psi_1 \qquad + \qquad \psi_2.$$

To solve the equation

$$\Delta \varphi_3 = \varphi_0 - \varphi_2 - \psi = \psi_1 + \psi_2,$$

we will solve separately the equations

$$\Delta \Psi_1 = \psi_1 \quad \text{and} \quad \Delta \Psi_2 = \psi_2,$$

and then add these solutions together to get

$$\varphi_3 = \Psi_1 + \Psi_2.$$

**First piece.**

If we define
$$\Psi_1(\alpha) = (1/4\pi)\, \alpha / (e^{2\alpha} - 1) - (1/4\pi)\, \text{csch}\, \alpha,$$

then
$$\Delta \Psi_1 = (1/4\pi) - (1/4\pi) \sinh \alpha / (1 + \cosh \alpha) = \psi_1.$$

**Second piece.**

We want to solve the equation
$$\Delta \Psi_2 = (1/4\pi)\, \alpha / (1 + \cosh \alpha) = \psi_2.$$

Following the procedure used on $S^3$ in section 31, we find that
$$\Psi_2 = (1/4\pi)\, \text{csch}\, \alpha + (1/4\pi) \int_0^\alpha \left((\alpha / \sinh \alpha) - (\alpha^2 / 2 \sinh^2 \alpha)\right) d\alpha.$$

The integral in the second term is non-elementary and involves dilogarithms, so we simply leave it as is.

Putting this all together, we have that
$$\varphi_3 = \Psi_1 + \Psi_2$$
$$= (1/4\pi)\, \alpha / (e^{2\alpha} - 1) - (1/4\pi)\, \text{csch}\, \alpha$$
$$\quad + (1/4\pi)\, \text{csch}\, \alpha + (1/4\pi) \int_0^\alpha \left((\alpha / \sinh \alpha) - (\alpha^2 / 2 \sinh^2 \alpha)\right) d\alpha$$
$$= (1/4\pi)\, \alpha / (e^{2\alpha} - 1) + (1/4\pi) \int_0^\alpha \left((\alpha / \sinh \alpha) - (\alpha^2 / 2 \sinh^2 \alpha)\right) d\alpha,$$

thanks to one cancellation.

This completes the proof of Theorem 3, formula (3) for the Green's operator in $H^3$.

# 52. Non-uniqueness of the Green's operator in $H^3$.

We have been seeking a formula for the vector-valued Green's operator in $H^3$ of the type

$$Gr(V)(y) = A(V, \varphi_2)(y) + G(V, \varphi_3)(y)$$

$$= \int_{H^3} P_{yx}V(x)\, \varphi_2(x, y)\, dx + \nabla_y \int_{H^3} P_{yx}V(x) \bullet \nabla_y \varphi_3(x, y)\, dx .$$

The ambiguity in the choice of $\varphi_2(\alpha)$ appeared in section 49, and is of the form $(1/4\pi)\, C / (1 + \cosh \alpha)$. So let's set

$$\phi_2 = 1 / (1 + \cosh \alpha) ,$$

and think of this as the difference between two competing choices of $\varphi_2$. We will have

$$\nabla_y \times A(V, \phi_2)(y) = 0 ,$$

since the corresponding integrals for the two choices of $\varphi_2$ both have the same curl, namely $-BS(V)(y)$.

Let $\phi_3$ denote the difference between the corresponding choices of $\varphi_3$. This kernel function satisfies the equation

$$\Delta_y \{A(V, \phi_2)(y) + G(V, \phi_3)(y)\} = 0 .$$

We follow the routine of sections 50 and 51 and learn that

$$\Delta \phi_3 = 1/(1 + \cosh \alpha) .$$

We solve this for $\phi_3$ and get

$$\phi_3'(\alpha) = (\sinh \alpha - \alpha) / \sinh^2 \alpha ,$$

and then

$$\phi_3(\alpha) = 2\alpha / (e^{2\alpha} - 1) + 2 \ln(e^\alpha / (e^\alpha + 1)) .$$

Putting this together, we get an explicit formula for the difference $\delta Gr$ between two competing vector-valued Green's operators in $H^3$,

$$\delta Gr(V)(y) = A(V, \phi_2)(y) + G(V, \phi_3)(y)$$

$$= \int_{H^3} P_{yx}V(x) \, \phi_2(x, y) \, dx + \nabla_y \int_{H^3} P_{yx}V(x) \bullet \nabla_y \phi_3(x, y) \, dx$$

$$= \int_{H^3} P_{yx}V(x) \, \big(1 / (1 + \cosh \alpha)\big) \, dx$$

$$+ \nabla_y \int_{H^3} P_{yx}V(x) \bullet \big((\sinh \alpha - \alpha) / \sinh^2 \alpha\big) \, \nabla_y \alpha \, dx \, .$$

We certainly know that $\delta Gr(V)$ is a harmonic vector field, that is, its Laplacian $\Delta_y \, \delta Gr(V)(y) = 0$. But recall from the construction that the first integral in the definition of $\delta Gr(V)$ is curl-free. Since $H^3$ is simply connected, this integral is actually a gradient vector field. The second integral is visibly a gradient field, so $\delta Gr(V)$ is actually a gradient field. Since it is harmonic, it is the gradient of some harmonic function on $H^3$, that is,

$$\delta Gr(V)(y) = \nabla_y(\text{some harmonic function of } y) \, .$$

Further computation reveals that this harmonic function on $H^3$ is in general *not* constant, and hence the arbitrary choice of the constant C in the definition of $\varphi_2$ in section 49 *does* have an effect on $Gr(V)$. As a consequence, the vector-valued Green's operator in $H^3$, *even restricted to the form that we are seeking,* is not unique.

# E. ELECTRODYNAMICS AND THE GAUSS LINKING INTEGRAL ON $S^3$ AND $H^3$

## 53. Proof of Theorem 4.

Recall from section 7 our assertion in Theorem 4 that the classical Maxwell equations,

$$\nabla \bullet E = \rho \qquad\qquad \nabla \times E = 0$$

$$\nabla \bullet B = 0 \qquad\qquad \nabla \times B = V + \partial E/\partial t,$$

hold on $S^3$ and in $H^3$, just as they do in $R^3$.

We have

$$\nabla_y \bullet E(\rho)(y) = \nabla_y \bullet \nabla_y \int \rho(x)\, \varphi_0(x, y)\, dx$$

$$= \Delta_y \int \rho(x)\, \varphi_0(x, y)\, dx$$

$$= \rho(y),$$

because $\varphi_0$ is the fundamental solution of the scalar Laplacian.

We have $\nabla \times E = 0$ because $E$ is a gradient field.

The equations

$$\nabla \bullet B = 0 \quad\text{and}\quad \nabla \times B = V + \partial E/\partial t$$

were two of the properties required of $B = BS(V)$ in section 4. They were proved on $S^3$ in left-translation format in section 12, on $S^3$ in parallel transport format in sections 23 and 24, and in $H^3$ in parallel transport format in sections 44 and 45.

This completes the proof of Theorem 4.

## 54. Proof scheme for Theorem 1.

Recall from section 2 what we imagine to be Gauss's second line of reasoning for obtaining his linking integral. *Run a current through the first loop, and calculate the circulation of the resulting magnetic field around the second loop. By Ampere's Law, this circulation is equal to the total current enclosed by the second loop, which means the current flowing along the first loop multiplied by the linking number of the two loops. Then the Biot-Savart formula for the magnetic field leads directly to Gauss's linking integral.*

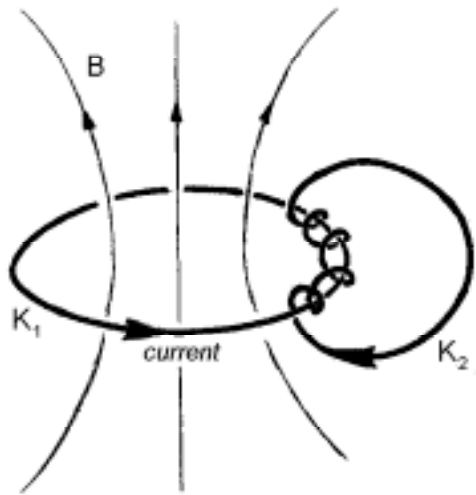

**Linking numbers via Ampere's Law**

The same line of reasoning in $S^3$ and $H^3$ will lead us in the next two sections from the three formulas for the magnetic field given in Theorem 2 to the three linking integrals given in Theorem 1.

But there is an issue which must be dealt with first. We have proved Theorem 2 only when V is a smooth, compactly supported current distribution on $R^3$, $S^3$ or $H^3$. As things stand, when currents run through wires, we have not proved the corresponding formula for the resulting magnetic field, nor the accompanying Ampere's Law.

We will get around this by observing that a field defined for a current running in a wire loop is the limit, in an appropriate sense, of magnetic fields caused by currents running in thin tubular neighborhoods of the loop. Then the Ampere's Law that we want will simply be the limit of Ampere's Laws that we already know, as the thickness of the neighborhoods goes to zero.

The argument is elementary, so we just outline it and leave the details to the reader.

Let K be a simple closed curve in an oriented Riemannian 3-manifold M, which in the present case is either $R^3$, $S^3$ or $H^3$. Let $r_0 > 0$ be a small real number, and let $N_{r_0}$ be the $r_0$-neighborhood of K in M, that is, the set of all points in M whose distance from K is less than $r_0$. We choose $r_0$ small enough so that $N_{r_0}$ is a traditional embedded tubular neighborhood of K.

Let K = {X(s), 0 ≤ s ≤ L} be a parametrization of K by arc length, and let T(s) = dX/ds = X'(s) be the resulting unit tangent vector field along K. Extend this vector field to an orthonormal frame field T(s), N(s), B(s) along K. This is possible because the normal bundle of a smooth curve in an oriented manifold is always trivial. The notation is chosen to remind us of the Frenet frame, consisting of the tangent, principal normal and binormal vectors along a smooth curve of nowhere vanishing curvature in Euclidean 3-space $R^3$.

Given this orthonormal frame field along K, we can use the exponential map of the Riemannian manifold M to parametrize the tubular neighborhood $N_{r_0}$ by the product of K and a Euclidean disk of radius $r_0$ as follows:

$$(s, u, v) \to X(s, u, v) = \exp_{X(s)}(u\, N(s) + v\, B(s)),$$

where $0 \leq s \leq L$ and $0 \leq u^2 + v^2 < r_0^2$. In $R^3$, this would simplify to

$$X(s, u, v) = X(s) + u\, N(s) + v\, B(s).$$

Now let $V(s, u, v) = \partial X(s, u, v)/\partial s$ be the tangent vector field to the s-coordinate curves, each of which is a simple closed curve having constant u and v coordinates. Since the exponential map takes rays in the normal plane to K at X(s) to geodesics orthogonal to K, each of the s-coordinate curves remains at the constant geodesic distance $(u^2 + v^2)^{1/2}$ from K.

**LEMMA.** The vector field $V(s, u, v) = \partial X(s, u, v)/\partial s$ can be rescaled to make it divergence-free.

The proof is elementary and we leave it to the reader.

It is easy to further rescale and concentrate the flow $V$ in $N_{r_0}$ to a flow $f_\varepsilon V$ that is supported in thinner and thinner tubular neighborhoods $N_\varepsilon$ of $K$, yet always has unit flux through the cross-sectional disks.

As $\varepsilon \to 0$, we regard the fields $f_\varepsilon V$ supported inside the neighborhoods $N_\varepsilon$ of $K$ as closer and closer approximations to the unit tangent field $X'(s) = T(s)$ along $K$.

We now change notation slightly for convenience.

Let $K = \{x(s)\}$ be a smooth simple closed curve in $R^3$, $S^3$ or $H^3$, and let $V = dx/ds$ be the unit tangent vector field along $K$. Let $BS(V)$ denote the field defined in the complement of $K$ by using any of the Biot-Savart formulas in $R^3$, $S^3$ or $H^3$ and by carrying out the required integrations with respect to arc length along $K$, rather than with respect to volume in the ambient space.

Let $V_\varepsilon$ denote a divergence-free field supported inside the tubular neighborhood $N_\varepsilon$ of $K$ and having unit flux through the cross-sectional disks, as constructed above. In particular, all the vector fields $V_\varepsilon$ are rescalings of a single vector field. Let $BS(V_\varepsilon)$ denote the corresponding magnetic field defined throughout the ambient space by using the same Biot-Savart formula as for $V$.

**LEMMA.** As $\varepsilon \to 0$, the magnetic fields $BS(V_\varepsilon)$ converge uniformly to the field $BS(V)$ on any compact subset of the complement of $K$.

Again, the proof of this is straightforward, and we leave the details to the reader.

As a consequence, Ampere's law for the field $BS(V)$ is the limit of Ampere's Laws for the fields $BS(V_\varepsilon)$ as $\varepsilon \to 0$. In particular, it can and will be applied to a loop in the complement of $K$ which bounds a surface cutting through $K$.

## 55. Proof of Theorem 1, formula (1).

Let $K_1 = \{x(s)\}$ and $K_2 = \{y(t)\}$ be disjoint oriented smooth closed curves in the unit 3-sphere $S^3$. We must show that in left-translation format, we have

$$Lk(K_1, K_2) = 1/(4\pi^2) \int_{K_1 \times K_2} L_{yx^{-1}} \, dx/ds \times dy/dt \bullet \nabla_y \, \varphi(x, y) \, ds \, dt$$

$$- 1/(4\pi^2) \int_{K_1 \times K_2} L_{yx^{-1}} \, dx/ds \bullet dy/dt \, ds \, dt \,,$$

where $\varphi(\alpha) = (\pi - \alpha) \cot \alpha$.

For the moment, we assume that the curve $K_1$ is simple, and that its parametrization $x = x(s)$ is by arc-length, so that $V(x) = dx/ds$ is a unit tangent vector field along $K_1$.

Then we take Theorem 2, formula (1) for the Biot-Savart operator in left-translation format and, given the discussion in the preceding section, apply it to the current flow $V$ along the loop $K_1$ to get

$$BS(V)(y) = \int_{K_1} L_{yx^{-1}} V(x) \times \nabla_y \, \varphi_0(x, y) \, ds$$

$$- 1/4\pi^2 \int_{K_1} L_{yx^{-1}} V(x) \, ds$$

$$+ 2 \nabla_y \int_{K_1} L_{yx^{-1}} V(x) \bullet \nabla_y \varphi_1(x, y) \, ds \,,$$

where $\varphi_0(\alpha) = (-1/4\pi^2)(\pi - \alpha) \cot \alpha$ and $\varphi_1(\alpha) = (-1/16\pi^2) \, \alpha \, (2\pi - \alpha)$.

Note that

$$L_{yx^{-1}} V(x) \bullet \nabla_y \varphi_1(x, y) = - V(x) \bullet \nabla_x \varphi_1(x, y)$$

$$= - dx/ds \bullet \nabla_x \varphi_1(x(s), y)$$

$$= - (d/ds) \, \varphi_1(x(s), y) \,,$$

whose integral around the loop $K_1$ is zero. So we can drop the third integral and are left with

$$BS(V)(y) = \int_{K_1} L_{yx^{-1}} V(x) \times \nabla_y \varphi_0(x, y) \, ds$$

$$- 1/4\pi^2 \int_{K_1} L_{yx^{-1}} V(x) \, ds \,.$$

Then we compute the circulation of $BS(V)$ around the loop $K_2$, as follows.

$\int_{K_2} BS(V)(y) \bullet dy/dt \, dt$

$= \int_{K_2} \left( \int_{K_1} L_{yx^{-1}} V(x) \times \nabla_y \varphi_0(x, y) \, ds - 1/4\pi^2 \int_{K_1} L_{yx^{-1}} V(x) \, ds \right) \bullet dy/dt \, dt$

$= \int_{K_1 \times K_2} L_{yx^{-1}} dx/ds \times \nabla_y \varphi_0(x, y) \bullet dy/dt \, ds \, dt$

$\qquad - 1/4\pi^2 \int_{K_1 \times K_2} L_{yx^{-1}} dx/ds \bullet dy/dt \, ds \, dt$

$= -\int_{K_1 \times K_2} L_{yx^{-1}} dx/ds \times dy/dt \bullet \nabla_y \varphi_0(x, y) \, ds \, dt$

$\qquad - 1/4\pi^2 \int_{K_1 \times K_2} L_{yx^{-1}} dx/ds \bullet dy/dt \, ds \, dt \, ,$

where we still have $\varphi_0(\alpha) = (-1/4\pi^2)(\pi - \alpha) \cot \alpha$,

$= 1/4\pi^2 \int_{K_1 \times K_2} L_{yx^{-1}} dx/ds \times dy/dt \bullet \nabla_y \varphi_0(x, y) \, ds \, dt$

$\qquad - 1/4\pi^2 \int_{K_1 \times K_2} L_{yx^{-1}} dx/ds \bullet dy/dt \, ds \, dt \, ,$

where now $\varphi_0(\alpha) = (\pi - \alpha) \cot \alpha$.

We apply Ampere's Law to the field $BS(V)$, as indicated in the previous section, and conclude that

$Lk(K_1, K_2) = 1/4\pi^2 \int_{K_1 \times K_2} L_{yx^{-1}} dx/ds \times dy/dt \bullet \nabla_y \varphi_0(x, y) \, ds \, dt$

$\qquad - 1/4\pi^2 \int_{K_1 \times K_2} L_{yx^{-1}} dx/ds \bullet dy/dt \, ds \, dt \, .$

The assumption that the loop $K_1$ is simple, which we needed to get an embedded tubular neighborhood, can now be dropped, since the linking number does not change as either curve crosses itself, and the value of the integral depends continuously on the curves $K_1$ and $K_2$.

The assumption that the loop $K_1$ is parametrized by arc length can also be dropped, since the two integrals on the right hand side above are independent of parametrization of either curve, as long as these parametrizations respect the give orientations.

This completes the proof of Theorem 1, formula (1) for the Gauss linking integral on $S^3$ in left-translation format.

## 56. Proof of Theorem 1, formulas (2) and (3).

We start on $S^3$, and begin with the same set-up and notation used in the preceding section. We must show that in parallel transport format, we have

$$\text{Lk}(K_1, K_2) = 1/4\pi^2 \int_{K_1 \times K_2} P_{yx} \, dx/ds \times dy/dt \bullet \nabla_y \varphi(x, y) \, ds \, dt \, ,$$

where $\varphi(\alpha) = (\pi - \alpha) \csc \alpha$.

We take Theorem 2, formula (2) for the Biot-Savart operator in parallel transport format and apply it to the current flow $V$ along the loop $K_1$ to get

$$\text{BS}(V)(y) = \int_{K_1} P_{yx} V(x) \times \nabla_y \varphi(x, y) \, ds \, ,$$

where $\varphi(\alpha) = (-1/4\pi^2)(\pi - \alpha) \csc \alpha$.

Then we compute the circulation of $\text{BS}(V)$ around the loop $K_2$, as follows.

$\int_{K_2} \text{BS}(V)(y) \bullet dy/dt \, dt =$

$$= \int_{K_2} \left( \int_{K_1} P_{yx} \, dx/ds \times \nabla_y \varphi(x, y) \, ds \right) \bullet dy/dt \, dt$$

$$= \int_{K_1 \times K_2} P_{yx} \, dx/ds \times \nabla_y \varphi(x, y) \bullet dy/dt \, ds \, dt$$

$$= -\int_{K_1 \times K_2} P_{yx} \, dx/ds \times dy/dt \bullet \nabla_y \varphi(x, y) \, ds \, dt \, ,$$

where we still have $\varphi(\alpha) = (-1/4\pi^2)(\pi - \alpha) \csc \alpha$,

$$= 1/4\pi^2 \int_{K_1 \times K_2} P_{yx} \, dx/ds \times dy/dt \bullet \nabla_y \varphi(x, y) \, ds \, dt \, ,$$

where now $\varphi(\alpha) = (\pi - \alpha) \csc \alpha$.

As before, we apply Ampere's Law to the field $\text{BS}(V)$ and conclude that

$$\text{Lk}(K_1, K_2) = 1/4\pi^2 \int_{K_1 \times K_2} P_{yx} \, dx/ds \times dy/dt \bullet \nabla_y \varphi(x, y) \, ds \, dt \, .$$

Then, as in the preceding section, we can drop the assumptions that the loop $K_1$ is simple and is parametrized by arc length, completing the proof of Theorem 1, formula (2) for the Gauss linking integral on $S^3$ in parallel transport format

The proof of Theorem 1, formula (3) for the Gauss linking integral on $H^3$ in parallel transport format is the same, and we are done.

University of Pennsylvania
Philadelphia, PA  19104

deturck@math.upenn.edu
gluck@math.upenn.edu